\documentclass[11pt, a4paper]{amsart}
\pdfoutput=1

\usepackage{amsmath, amssymb, amsthm, bbm, multicol, rotating, tabularx, enumitem}
\usepackage[all]{xy}
\usepackage[a4paper, left=2.5cm, right=3cm, top=3cm, bottom=3cm]{geometry}
\usepackage[english]{babel}

\usepackage{graphicx, color}
\usepackage{ucs}
\usepackage[utf8x]{inputenc}
\usepackage[T1]{fontenc}

\usepackage[pdfpagelabels]{hyperref}

\theoremstyle{plain}
\newtheorem{thm}{Theorem}[section]
\newtheorem*{thm*}{Theorem}
\newtheorem{prop}[thm]{Proposition}
\newtheorem*{prop*}{Proposition}
\newtheorem{lemma}[thm]{Lemma}
\newtheorem*{lemma*}{Lemma}
\newtheorem{corollary}[thm]{Corollary}
\newtheorem*{corollary*}{Corollary}
\newtheorem{property}[thm]{}
\newtheorem{conj}[thm]{Conjecture}

\newtheorem{defthm}[thm]{Theorem and Definition}

\newtheorem*{ThmConv*}{Theorem \ref{Thm_convexityGeneral}}
\newtheorem*{ThmBNpair*}{Theorem \ref{Thm_BNpair}}

\theoremstyle{definition}
\newtheorem{definition}[thm]{Definition}
\newtheorem*{definition*}{Definition}
\newtheorem{example}[thm]{Example}

\theoremstyle{remark}
\newtheorem{remark}[thm]{Remark}
\newtheorem{notation}[thm]{Notation}

\newcommand{\R}{\mathbb{R}} 

\newcommand{\Z}{\mathbb{Z}}
\newcommand{\Q}{\mathbb{Q}}
\newcommand{\define}{\mathrel{\mathop:}=}
\newcommand{\ddefine}{\mathrel{=\mathop:}}

\newcommand{\MS}{\mathbb{A}} 
\newcommand{\App}{\mathcal{A}}
\newcommand{\seg}{\mathrm{seg}} 
\newcommand{\binfinity}{\partial_\App} 

\newcommand{\RS}{\mathrm{R}}
\newcommand{\rk}{\mathrm{rank}}
\newcommand{\QQ}{{\bf\mathrm{Q}}} 
\newcommand{\fcw}[1][\alpha]{{\omega}_{#1}}
\newcommand{\sW}{\overline{W}} 
\newcommand{\aW}{W} 
\newcommand{\WT}{\overline{W}T}
\newcommand{\cf}{c_0} 
\newcommand{\Cf}{\mathcal{{C}}_{f}} 
\newcommand{\Cfm}{-\mathcal{{C}}_{f}} 
\newcommand{\Cp}{\mathcal{{C}}_{p}} 

\newcommand{\dconv}{\mathrm{{conv}}} 
\newcommand{\Hdual}{\overline{H}}
\newcommand{\Ghat}{\hat{\mathcal{G}}} 
\newcommand{\G}{\mathcal{G}} 

\newcommand{\one}{\mathbbm{1}} 

\newcommand{\lb}{\langle} 
\newcommand{\rb}{\rangle} 

\DeclareMathOperator{\Span}{span}

\numberwithin{equation}{thm}

\begin{document}

\hypersetup{pdfauthor={Petra Hitzelberger},pdftitle={Non-discrete affine buildings and Convexity}}.
\title{Non-discrete affine buildings and convexity}
\author{Petra Hitzelberger}
\thanks{{\it Address of the author:}\
Fachbereich Mathematik und Informatik, Universit\"at M\"unster,
Einsteinstrasse~62, 48149 M\"unster, Germany} 
\date{ June 2009 }
\maketitle

\begin{abstract}
Affine buildings are in a certain sense analogs of symmetric spaces. It is therefore natural to try to find analogs of results for symmetric spaces in the theory of buildings. In this paper we prove a version of Kostant's convexity theorem for thick non-discrete affine buildings. Kostant proves that the image of a certain orbit of a point $x$ in the symmetric space under a projection onto a maximal flat is the convex hull of the Weyl group orbit of $x$. We obtain the same result for a projection of a certain orbit of a point in an affine building to an apartment. The methods we use are mostly borrowed from metric geometry. Our proof makes no appeal to the automorphism group of the building. However the final result has an interesting application for groups acting nicely on non-discrete buildings, such as groups admitting a root datum with non-discrete valuation. Along the proofs we obtain that segments are contained in apartments and that certain retractions onto apartments are distance diminishing. 
\end{abstract}

%
%
%
%
\section{Introduction}
\label{sec_introduction}

Kostant's convexity theorem for symmetric spaces, proven in \cite{Kostant}, describes the image of a certain orbit under a projection on a maximal flat as a convex set. His result is a generalization of a well known theorem of Schur \cite{Schur}. The precise geometric statement is as follows:

Let $G\!/\!K$ be a symmetric space and $T$ a maximal flat of $G\!/\!K$. Then there is a natural action of a spherical Weyl group $\sW$ on $T$ with fixed point $0\in T$. Write $\pi$ for the Iwasawa projection of $G$ onto $T$.
Kostant proves that the image of the orbit $K.x$ of an element $x\in T$ under the Iwasawa projection is precisely the convex hull of the Weyl group orbit $\sW.x$.  In terms of groups his result provides a criterion for the non-emptiness of intersections of certain double cosets of group elements.

Let $G$ be a non-compact semi-simple Lie-group with Iwasawa decomposition $G=UTK$ with $U$  unipotent, $T$ abelian and $K$ compact. Geometrically  left-cosets of elements of $T$ correspond to points in a certain maximal flat of the symmetric space $G\!/\! K$. Kostant's theorem translates o the fact that intersections of double cosets of the form $KbK\cap Ub'K$ are non-empty if and only if $b'K$ is contained in the convex hull of the point $bK$. 

Affine buildings are in a certain sense analogs of symmetric spaces. Similarly symmetric spaces are important for the classification of semisimple Lie groups, so do affine buildings play a major role in the classification of semisimple algebraic groups defined over fields with valuation. Part of the analogy in terms of geometry is as follows: Maximal flats in symmetric spaces correspond to apartments in buildings. Both of them admit an action of a spherical Weyl group. 

The notion of an ``Iwasawa projection'' $\rho$ onto an apartment $A$ does make sense in a building, too. In terms of groups it is defined precisely in the same way as in the context of semi-simple Lie groups, but there is as well a definition using the geometry of an affine building. The orbit $K.x$ in Kostant's result corresponds, when talking about buildings, to the preimage of $\sW.x$ under a second type of retraction onto $A$, which we will denote by $r$. Hence we might ask again whether the projection of this set onto $A$ is a convex hull of the Weyl group orbit of $x$. Or, spoken in group language, whether for algebraic groups (to be precise groups with affine and split spherical BN-pair) the same criterion guarantees the intersections of double cosets to be non-empty.

In the simplicial case this question was answered in \cite{Convexity}. The purpose of the present paper is to prove a convexity result in the spirit of Kostant's for a class of spaces more general than simplicial affine buildings.

\subsection*{Generalized affine buildings}
Simplicial affine buildings, which are a subclass of the geometric objects studied in this paper, were introduced by Bruhat and Tits in \cite{BruhatTits} as spaces associated to semisimple algebraic groups defined over fields with discrete valuations. 

In \cite{TitsComo} and \cite{BruhatTits, BruhatTits2} they were generalized allowing fields with non-discrete (non-Archimedean) valuations rather than discrete ones. The arising geometries do no longer carry a simplicial structure and are nowadays usually called \emph{non-discrete affine buildings} or \emph{$\R$-buildings}. In \cite{TitsComo}, $\R$-buildings were axiomatized and for sufficiently large rank classified under the name \emph{syst\`eme d'appartements}. 

Finally, in \cite{Bennett, BennettDiss} Bennett introduced a class of $\Lambda$-metric spaces called \emph{affine $\Lambda$-buildings} using axioms similar to the ones in \cite{TitsComo}. Examples of these spaces arise from simple algebraic groups defined over fields with valuations taking their values in an arbitrary ordered abelian group $\Lambda$. Bennett was able to prove that affine $\Lambda$-buildings again have simplicial spherical buildings at infinity and made major steps towards their classification. 
To be more precise a generalized affine building is a set $X$ together with a collection of maps $\App$ called \emph{atlas}. Each $f\in\App$ is an injective map from a (fixed) model apartment $\MS$ to $X$. The images $f(\MS)$ are called \emph{apartments} of $X$. As a set $X$ is the union of its apartments, which need to satisfy certain axioms in addition. Compare Definition~\ref{Def_LambdaBuilding}.

For generalized affine buildings one can again define an ``Iwasawa projection'' $\rho: X \to A$ onto an apartment $A$ and one can, as already mentioned above, define a second type of retraction $r: X \to A$ whose preimage of $\sW.x$ corresponds precisely to the $K$-orbit of $x$. Therefore it is natural to ask whether a a convexity theorem exists for this more general class of affine buildings.

Throughout this text we will refer to affine $\Lambda$-buildings as \emph{generalized affine buildings} to avoid the appearance of the group $\Lambda$ in the name.
Note that the class of generalized affine buildings does not only include all previously known classes of (non-discrete) affine buildings, but does also generalize leafless \emph{$\Lambda$-trees}. These trees are simply the generalized affine buildings of dimension one. 

\subsection*{Convexity} 
Let $(X,\App)$ be a thick generalized affine building, as defined in \ref{Def_LambdaBuilding} and \ref{Def_thick} with model space $\MS$. There is an action of an affine Weyl group $\aW$ on $\MS$. The stabilizer of $0\in\MS$ in $\aW$ can naturally be identified with the spherical Weyl group $\sW$. We fix a fundamental domain of the action of $\sW$ on $\MS$ and call it the \emph{fundamental Weyl chamber} $\Cf$. Weyl chambers in $\MS$ are images of $\Cf$ under the affine Weyl group $\aW$ and Weyl chambers in $X$ are images of Weyl chambers in $\MS$. Fixing a chart $f$ of an apartment $A$ in $X$ it therefor makes sense to talk about an origin $0$ and a fundamental Weyl chamber $\Cf$ in $A$. We say that two Weyl chambers based at the same vertex are \emph{equivalent} if they intersect in a set with nonempty relative interior. The equivalence class of a Weyl chamber $S$ based at $x$ is called the \emph{germ} of $S$ at $x$.

To state the main theorem we need to introduce two retractions onto the given apartment $A$.
The first one, denoted by $r$, is defined with respect to the germ of the fundamental Weyl chamber $\Cf$ in $A$.  It preserves distances to $0$ and its restriction to apartments containing the germ of $\Cf$ at $0$ is an isomorphism onto $A$. The inverse image $r^{-1}(\sW.x)$ of the Weyl group orbit of $x$ under $r$ corresponds precisely to the orbit $K.x$ in Kostant's setting. The second retraction $\rho : X\mapsto A$ is sometimes (mostly when talking about algebraic groups) called ``Iwasawa projection'' onto $A$. The geometric definition of $\rho$ is given with respect to the parallel class $\partial (\Cfm)$ of the opposite of the fundamental Weyl chamber.  Here two Weyl chambers are \emph{parallel} if their intersection contains a Weyl chamber. We demand that the restriction of $\rho$ to an apartment containing a Weyl chamber parallel to $\Cfm$ is an isomorphism onto $A$. As it turns out this leads to a well defined retraction of $X$ onto $A$. We obtain the following theorem. 

\begin{ThmConv*} 
Given a vertex $x$ in $A$ one has
$$
\rho(r^{-1}(\sW.x)) =\dconv(\sW.x).
$$
\end{ThmConv*}

Restating this in terms of a group $G$ acting ``nicely'' on a thick generalized affine building $X$, one obtains the following result about non-emptiness of intersections of double cosets in $G$.
Denote by $K$ the stabilizer of the origin $0$ in $G$ and assume that it is transitive on the apartments containing $0$. Let further $B=G_c$ be the stabilizer of the equivalence class $c=\partial( \Cfm)$ and assume that $B$ splits as $B=UT$, where $T$ is the group of translations in $A$ and $U$ acts simply transitive on the apartments containing $c$ at infinity. Then 

\begin{ThmBNpair*}
For all $tK\in A$ we have 
$$
\rho(KtK)=\dconv(\sW.tK)
$$
or, since $\rho^{-1}(t'K)=Ut'K$, equivalently
$$
Ut'K \cap KtK \neq \emptyset  \:\Longleftrightarrow\: t'K \in \dconv(\sW.tK).
$$
\end{ThmBNpair*}

Techniques used to prove Theorem \ref{Thm_convexityGeneral} are geometric properties of generalized affine buildings and methods borrowed from metric geometry. An important  idea is inspired by and, with enough technical effort, adapted from a result of Parkinson and Ram proven in \cite{ParkinsonRam}. They give a combinatorial proof of the existence of certain positively folded galleries. The main idea of their proof can be modified to obtain a result on retractions onto apartments of generalized affine buildings.

\subsection*{Outline of proof}
The first problem arising concerns the two retractions $r$ and $\rho$. 
In order to be able to define them in the general setting of the present paper, we first have to establish several structural results on the local and global behavior of generalized affine buildings. This is done in Section~\ref{Sec_local-global}. The definition of the retractions can be found in Section~\ref{Sec_retractions} 
In this context we introduce residues, which are as sets simply the collection of all germs based at the same vertex, but carry the structure of spherical buildings. One can think of residues as ``tangent spaces'' at points in $X$.

Further we want to prove that both retractions do not increase distances between arbitrary points in the building. This fact is much easier to prove in the simplicial case. In order to verify that they  are distance-non-increasing (see Corollary \ref{Cor_rho-distancediminishing}), we need to generalize Lemma 7.4.21 of Bruhat and Tits \cite{BruhatTits} which is a covering property of segments in buildings. The proof of Lemma 7.4.21 given in \cite{BruhatTits} uses compactness arguments of $\R$-metric spaces which cannot be applied in our setting. In Section~\ref{Sec_FC} we prove these properties for generalized affine buildings.

We are then ready to prove Theorem \ref{Thm_convexityGeneral}. The major problem occurs in the proof of the fact that every element of the convex set $\dconv(\sW.x)$ has a preimage under $\rho$ which is contained in the set $r^{-1}(\sW.x)$. 
As already mentioned earlier this is done by modifying an idea of Parkinson and Ram \cite{ParkinsonRam}. Given an element $y$ of $\dconv(\sW.x)$ we define, in the proof of Proposition~\ref{Prop_preimage}, a sequence of points $y_i$ depending on a chosen presentation of the longest element $w_0$ of the spherical Weyl group $\sW$. This sequence of points helps to define more or less explicitly a preimage of $y$ under $\rho$, which is by construction contained in $r^{-1}(\sW.x)$. As pretty as the main idea may be as technical is the actual proof. For the convenience of the reader we therefore repeat the underlying ideas without proof in Section~\ref{Sec_ParkinsonRam}.

To finish the proof of \ref{Thm_convexityGeneral} it remains to show that the image of the set $r^{-1}(\sW.x)$ is in fact contained in the convex hull of the Weyl group orbit of $x$. This is done in Proposition~\ref{Prop_image}. Methods used in the proof are borrowed from metric geometry and mimic differentiation. The ideas come from the similarity between germs of Weyl chambers and tangent vectors of curves in manifolds or metric spaces.

We have not mentioned so far that the notion of convexity used in the present paper is not the usual one where convex sets are defined to be finite intersections of half-apartments. However our notion of convexity, as defined in Definition~\ref{Def_GammaConvex}, generalizes the metric convex hull in terms of the Euclidean metric defined in the geometric realization of simplicial buildings. An apartment of an $\R$-building is naturally equipped with a Euclidean metric. The metric convex hull of a Weyl group orbit $\sW.x$ in such an apartment defined with respect to the Euclidean distance corresponds precisely to the convex hull as it is defined in \ref{Def_GammaConvex} (in case the building is equipped with the full affine Weyl group). However dealing with $\Lambda$-metric spaces there is nothing like a Euclidean distance. We therefore define a metric on apartments differently. Using this metric, described in definition \ref{Def_distance}, it is no longer true that the metric convex hull of $W.x$ equals the convex hull in the sense of Definition \ref{Def_GammaConvex}. Only the weaker observation \ref{Prop_segment} remains.

\subsection*{The paper is organized as follows} 
In section \ref{Sec_modelSpace} the building block of a generalized affine building, the so-called model space $\MS$ of an apartment, is defined. 
Generalized affine buildings are then defined in Section \ref{Sec_generalizedAffineBuildings}, where we also describe their local and global structure and prove preliminary results which are necessary for the definition of retractions. 

These retractions are then defined in Section \ref{Sec_retractions}. 
Local covering properties (generalizing a Lemma by Bruhat and Tits), which are used to prove that the retractions of the previous chapter are distance diminishing, are investigated in Section \ref{Sec_FC}. 

The following Section \ref{Sec_ParkinsonRam} might be skipped. Here we recall the convexity theorem proven in the setting of simplicial affine buildings and explain a geometrical construction of certain positively folded paths. This is done in order to make the technical proof of Theorem \ref{Thm_convexityGeneral}, given in Section \ref{Sec_convexityThm} and relying on the ideas of Section \ref{Sec_ParkinsonRam}, more approachable. 

An application to groups acting nicely enough on thick affine buildings is then given in Section \ref{Sec_application}. This application is similar to the one obtained in the simplicial case and the direct analog of the result by Kostant on non-emptiness of intersections of double cosets. 

Finally an open problem is discussed in Section \ref{Sec_looseEnds}.

\subsection*{Acknowledgments}

The author would like to thank Linus Kramer for many helpful discussions
and encouragement. We also thank James Parkinson for the reference to \cite{ParkinsonRam}. The author was partially supported by the \emph{Studienstiftung des deutschen Volkes} and the \emph{SFB 478 "Geometrische Strukturen in der Mathematik"} while working on this topic. This work is part of the author’s doctoral thesis at the Universität Münster

%
%
%
%

\section{The model space}\label{Sec_modelSpace}

Geometric realizations of simplicial affine buildings are metric spaces ``covered by'' Coxeter complexes which are isomorphic to a tiled $\R^n$. The basic idea of the generalization is to substitute the real numbers by a totally ordered abelian group $\Lambda$.

\subsection{Definition and basic properties}
\begin{definition}\label{Def_modelSpace}
\index{model space}
Let $\RS$ be a (not necessarily crystallographic) spherical root system and $F$ a subfield of $\R$ containing the set $\{\lb \alpha, \beta^\vee\rb : \alpha, \beta \in \RS\}$ of co-roots $\alpha^\vee$ evaluated on roots $\beta$. Assume that $\Lambda$ is a totally ordered abelian group admitting an $F$-module structure. The space 
$$
\MS(\RS,\Lambda) = \mathrm{span}_F(\RS)\otimes_F \Lambda
$$ 
is the \emph{model space} of a generalized affine building of type $\RS$. 
\end{definition}

We omit $F$ in the notation, since we can always choose $F$ to be $\Q(\{\lb \beta, \alpha^\vee\rb : \alpha, \beta \in \RS\})$. If $\RS$ is crystallographic then $F=\Q$ is  always a valid choice. If there is no doubt which root system $\RS$ and which $\Lambda$ we are referring to, we will abbreviate $\MS(\RS,\Lambda)$ by $\MS$.

\begin{remark}\label{Rem_coordinates}
\index{coordinates}
A fixed basis $B$ of the root system $\RS$ provides natural coordinates for the model space $\MS$. The space of formal sums 
$$
\left\{\sum_{\alpha\in B} \lambda_\alpha\alpha : \lambda_\alpha\in\Lambda\right\}
$$ 
is canonically isomorphic to $\MS$ and the evaluation $\lb\alpha, \beta^\vee\rb$ of co-roots $\beta^\vee$ on roots $\alpha$ can be extended linearly to its elements.
\end{remark}

\begin{definition}\label{Def_reflections}
\index{hyperplane}
An action of the spherical Weyl group $\sW$ on $\MS$ is defined as follows.
Let $\alpha, \beta \in\RS$, $c\in F$ and $\lambda\in\Lambda$, let $r_\alpha:\MS\mapsto\MS$ be the linear extension of 
$$r_\alpha(c\beta\otimes \lambda ):= c\, s_\alpha(\beta) \otimes \lambda$$
to $\MS$, where the \emph{reflection} $s_\alpha$ is defined by
\begin{equation}
s_\alpha(\beta) = \beta - 2\frac{(\alpha, \beta)}{(\alpha,\alpha)}\alpha.
\end{equation}
We call the fixed point set of $r_\alpha$ a \emph{hyperplane} or \emph{wall} and we denoted it by $H_\alpha$ or $H_{\alpha,0}$, since $H_\alpha=\{x \in\MS: \lb x, \alpha^\vee\rb =0\}$.
\end{definition}

A basis $B$ of $\RS$ determines a set of positive roots $\RS^+\subset \RS$. The subset 
$$
\Cf\define \{x\in\MS : \lb x,\alpha^\vee\rb \geq 0 \text{ for all } \alpha\in\RS^+\}
$$
of $\MS$ is the \emph{fundamental Weyl chamber} with respect to $B$, and is denoted by $\Cf$. 

\begin{definition}\label{Def_affineWG}
Given a non-trivial group of translations $T$ of $\MS$ which is normalized by $\sW$ we define the \emph{affine Weyl group} with respect to $T$ to be the semi direct product $\WT = \sW\rtimes T$. In case $T=\MS$ we call it the \emph{full affine Weyl group} and write $\aW=\sW\rtimes \MS$.
Elements of $T$ can be identified with points in $\MS$ by assigning to $t\in T$ the image of the origin $0$ under $t$. Given $k\in \MS$ we write $t_k$ for the translation defined by $0\mapsto k$.
\end{definition}

The actions of $\sW$ and $T$ on $\MS$ induce an action of $\WT$, respectively $\aW$, on $\MS$. 

\begin{notation}\label{Not_modelSpace}
In order to emphasize the freedom of choice for the translation part of the affine Weyl group, the model space $\MS(\RS, \Lambda)$ with affine Weyl group $\WT$ is referred to as $\MS(\RS, \Lambda, T)$.
\end{notation}

\begin{definition}\label{Def_specialHyperplane}
\index{{reflection}!{affine}}
\index{hyperplane}
An element $w$ of $\WT$ which can be written as $t \circ r_\alpha$ for some nontrivial $t\in T$ and $\alpha\in \RS$ with $r_\alpha(t)=-t$ is called \emph{(affine) reflection}.
A \emph{hyperplane} $H_r$ in $\MS$ is the fixed point set of an affine reflection $r \in \aW$. It is called \emph{special with respect to $T$} if $r\in \WT$. 
\end{definition}

\begin{remark}
Note that for any affine reflection there exists $\alpha\in\RS$ and $\lambda\in\Lambda$ such that the reflection is given by the following formula
$$
r_{\alpha,\lambda}(x) 
	= r_\alpha(x) + \frac{2\lambda}{(\alpha,\alpha)} \alpha, 
	\text{ for all } x\in \MS.
$$
Further, easy calculations imply that
$$
r_{\alpha, \lambda}(x) 
	= r_\alpha\left( x-\frac{\lambda}{(\alpha,\alpha)}\alpha \right) + \frac{\lambda}{(\alpha,\alpha)}\alpha ,
	\text{ for all } x\in \MS.
$$
The fixed point set $H_{\alpha, \lambda}$ of $r_{\alpha,k}$ is given by
$$
H_{\alpha,\lambda}=\left\{x \in\MS: \frac{(\alpha,\alpha)}{2}\lb x,\alpha^\vee\rb =\lambda \right\}.
$$

As in the classical case each hyperplane defines a positive and a negative  \emph{half-apartment} 
$$H_{\alpha,k}^+=\{x \in\MS: \frac{(\alpha,\alpha)}{2}\lb x,\alpha^\vee\rb \geq \lambda \}\; \text{ and }\;
H_{\alpha,k}^-=\{x \in\MS: \frac{(\alpha,\alpha)}{2}\lb x,\alpha^\vee\rb \leq \lambda \}.$$
\end{remark}

\begin{definition}
A vertex $x\in\MS$ is called \emph{special} if for each $\alpha\in\RS^+$ there exists a special hyperplane parallel to $H_{\alpha,0}$ containing $x$. Hence $x$ is the intersection of the maximal possible number of special hyperplanes.
\end{definition}

Note that the translates of $0$ by $T$ are a subset of the set of special vertices.

A \emph{Weyl chamber} in $\MS$ is an image of the fundamental Weyl chamber under the full affine Weyl group $\aW$. If $\Lambda=\R$ Weyl chambers are simplicial cones in the usual sense. Therefore Weyl chambers and their faces are called \emph{Weyl simplices}. The faces of co-dimension one are referred to as \emph{panels}.

Note that a Weyl chamber $S$ contains exactly one vertex $x$ which is the intersection of all bounding hyperplanes of $S$. We call it \emph{base point} of $S$ and say $S$ is \emph{based at $x$}.

The following proposition is used to introduce a second type of coordinates on $\MS$.

\begin{prop}{ \cite[Prop. 2.1]{Bennett}}\label{Prop_aboveHyperplane}
Given $\alpha\in \RS$ and $x\in\MS$. Then there exist a unique $m_\alpha \in H_{\alpha,0}$ and a unique $x^\alpha\in\Lambda$ such that 
$$x= m_\alpha + x^\alpha \alpha.$$
The value of $x^\alpha$ is $\frac{1}{2} \lb x, \alpha^\vee\rb$. Furthermore $x\in H_{\alpha,(\alpha,\alpha)x^\alpha}$.
\end{prop}
\begin{proof}
Define $x^\alpha = \frac{1}{2} \lb x,\alpha^\vee\rb$ and consider $m_\alpha=x- x^\alpha\alpha$. Then 
$$\lb  m_\alpha,\alpha^\vee\rb = \lb x, \alpha^\vee\rb - x^\alpha \lb \alpha,\alpha^\vee\rb = 0 $$
and $m_\alpha$ is contained in $H_{\alpha,0}$. 
It remains to prove uniqueness. 
Let $y^\alpha$ and $n_\alpha$ be such that 
$m_\alpha +x^\alpha\alpha = n_\alpha + y^\alpha\alpha.$
Then $n_\alpha=r_\alpha(n_\alpha) = r_\alpha(x) - y^\alpha r_\alpha(\alpha)$ and
$m_\alpha=r_\alpha(m_\alpha) = r_\alpha(x) - x^\alpha r_\alpha(\alpha)$. Therefore we have
$$
r_\alpha(x)-x^\alpha\alpha = m_\alpha + x^\alpha\alpha = x = n_\alpha + y^\alpha\alpha = r_\alpha(x) - y^\alpha\alpha 
$$ 
and conclude that $x^\alpha=y^\alpha$ and $m_\alpha=y_\alpha$.
\end{proof}

\begin{corollary}\label{Cor_hyperplaneCoordinates}
\index{hyperplane coordinates}
Let $\RS$ be a root system of rank $n$ with basis $B$. A point $x\in\MS$ is uniquely determined by the values $\{x^\alpha\}_{\alpha\in B}$, which will be called \emph{hyperplane coordinates} of $\MS$ with respect to $B$. 
\end{corollary}

\subsection{The metric structure of \texorpdfstring{ $\MS(\Lambda, \RS)$ } {the model space}}

The remainder of this section is used to define a $\aW$-invariant $\Lambda$-valued metric on the model space $\MS=\MS(\Lambda, \RS, T)$ of a generalized affine building and to discuss its properties.

\begin{definition}\label{Def_metric}
Let $\Lambda$ be a totally ordered abelian group and let $X$ be a set. A metric on $X$ with values in $\Lambda$, short a \emph{$\Lambda$-valued metric}, is a map $d:X\times X \mapsto \Lambda$ such that for all $x,y,z\in X$ the following axioms are satisfied
\begin{enumerate}
\item $d(x,y)=0$ if and only if $x=y$
\item $d(x,y)=d(y,x)$ and
\item the triangle inequality $d(x,z)+d(z,y)\geq d(x,y)$ holds.
\end{enumerate}
The pair $(X,d)$ is a \emph{$\Lambda$-metric space}.
\end{definition}

\begin{definition}
An \emph{isometric embedding} of a $\Lambda$-metric space $(X,d)$ into another $(X',d')$ is a map $f:X\rightarrow X'$ such that for all $x$ and $y$ in $X$ one has  $d(x,y)=d'(f(x),f(y))$. Such a map is necessarily injective, but need not be onto. If it is onto we call it an \emph{isometry} or an isomorphism of $\Lambda$-metric spaces.
\end{definition}

\begin{definition}\label{Def_distance}
Let $\MS=\MS(\Lambda,\RS, T)$ be as in \ref{Not_modelSpace}. The \emph{distance} of points $x$ and $y$ in $\MS$ is given by
$$d(x,y)= \sum_{\alpha\in\RS^+} \vert \lb y-x, \alpha^\vee \rb \vert.$$
\end{definition}

If $y-x$ is contained in the fundamental Weyl chamber the distance  $d(x,y)$ equals $2\lb y-x, \rho^\vee\rb$, where $\rho^\vee$ is half the sum of the positive co-roots.

Let $B$ be a basis of $\RS$. Choosing $\varepsilon_\alpha$ equal to $\frac{2}{(\alpha,\alpha)}$ and identifying the co-roots $\alpha^\vee$ with $\frac{2\alpha}{(\alpha, \alpha)}$ in the definition of the metric in \cite{BennettDiss} we precisely obtain the metric defined in \ref{Def_distance}.

The latter generalizes the chamber distance (or length of a translation) in a Coxeter Complex. This distance is defined as follows: Let $v$ be a vertex in a Euclidean Coxeter complex. Let the length $l(t_v)$ of a translation $t_v$ be the number of hyperplanes crossed by a minimal gallery from $x$ to $t_v(x)=x+v$. This number is given by the formula  $\frac{1}{2}\sum_{\alpha\in\RS} \vert \lb y-x, \alpha^\vee \rb \vert $.
The fact that $d$ is the direct generalization of a combinatorial length function justified, at least in my opinion, to make a specific choice of the $\varepsilon_i$ appearing in the definition of the metric as written in \cite{BennettDiss}.

\begin{prop}
The distance $d:\MS\times\MS\mapsto \Lambda$ defined in~\ref{Def_distance} is a $\aW$-invariant $\Lambda$-valued metric on $\MS$.
\end{prop}
\begin{proof}
By definition $d(x,y)=d(y,x)$ and $d(x,y)=0$ if and only if $x=y$, since otherwise, by Corollary 2.2 in \cite{Bennett}, one of the terms 
$\vert \lb z, \alpha^\vee\rb \vert$ would be strictly positive. 
It remains to prove that $d(x,y)+d(y,z)\geq d(x,z)$:
\begin{align*}
d(x,z) 	&= \sum_{\alpha\in\RS^+} \vert \lb  z+ (y-y) - x , \alpha^\vee \rb\vert \\
	&\leq \sum_{\alpha\in\RS^+} \big( \vert \lb y-x,\alpha^\vee\rb \vert +  \vert \lb z-y, \alpha^\vee\rb \vert \big) \\
	&= d(x,y)+d(y,z).
\end{align*}
Hence $d$ is a metric. We prove $\aW$-invariance: 
Let $t_a: x\mapsto x+a$ be a translation in $\aW$. Then
\begin{align*}
d(t_a(x), t_a(y)) &= \sum_{\alpha\in\RS^+} \vert \lb y+a - (x+a), \alpha^\vee \rb \vert
 	= \sum_{\alpha\in\RS^+} \vert \lb y-x, \alpha^\vee \rb \vert
	= d(x,y).
\end{align*}
Therefore $d$ is translation invariant. With $w\in \sW$ we have
\begin{align*}
d(w.x, w.y)
	& = \sum_{\alpha\in\RS^+} \vert \lb w.y-w.x,\alpha^\vee \rb \vert 
	    =  \frac{1}{2}\sum_{\alpha\in\RS} \vert \lb w.y-w.x,\alpha^\vee \rb \vert  \\
	& = \frac{1}{2}\sum_{\alpha\in\RS} \vert \lb y-x, (w^{-1}.\alpha)^\vee \rb \vert 
	    = \frac{1}{2}\sum_{\alpha\in\RS} \vert \lb y-x, \alpha^\vee \rb\vert \\
	& =d(x,y).
\end{align*}
The second last equation holds since $\sW$ permutes the roots in $\RS$.
Therefore $d$ is $\sW$ invariant and $\aW$-invariance follows.
\end{proof}

Recall the definition of the hyperplane coordinates $\{x^\alpha\}_{\alpha\in B}$ of $x$ with respect to a basis $B$ introduced in Corollary~\ref{Cor_hyperplaneCoordinates}.

\begin{prop}\label{Prop_distance0x}
Fix a basis $B$ of $\RS$ and let $x$ be an element of $\MS$. The distance $d(0,x)$ is uniquely determined by the hyperplane-coordinates $\{x^\alpha\}_{\alpha\in B}$ of $x$. 
With  $\alpha = \sum_{\beta\in B} p_\beta^\alpha\beta$ we have
$$
d(0,x)=\frac{1}{2} \sum_{\alpha \in \RS^+} \sum_{\beta\in B}  p_\beta^\alpha \; \vert x^\beta\vert.
$$
\end{prop}
\begin{proof}
Assume first that $x\in \Cf$. Then $x$ has hyperplane coordinates $\{x^\alpha\}_{\alpha\in B}$ defined in Corollary~\ref{Cor_hyperplaneCoordinates}, with $x^\alpha=\frac{1}{2}\lb x, \alpha^\vee\rb \geq 0$ for all $\alpha\in \RS^+$.  Hence, using $\alpha = \sum_{\beta\in B} p_\beta^\alpha\beta$, we have
\begin{align*}
d(x,0) 	&= \sum_{\alpha\in\RS^+} \vert \lb x, \alpha^\vee \rb\vert 
	 = \sum_{\alpha\in\RS^+} \sum_{\beta\in B} p_\beta^\alpha \,\lb x, \beta^\vee \rb
	 = \frac{1}{2} \sum_{\alpha\in\RS^+} \sum_{\beta\in B} p_\beta^\alpha  \,x^\beta.
\end{align*}
If $x$ is not $\Cf$ then $\aW$-invariance of $d$ implies that $d(0,x)=d(0,\overline{x})$, where $\overline{x} = w.x$ is the unique element of $\sW.x$ contained in $\Cf$ and $\overline{\beta}=w.\beta$. Further $\vert x^\beta \vert = \overline{x}^{\overline{\beta}}$ and  the assertion follows.
\end{proof}

\subsection{Convexity and parallelism}

As in the classical case, one can define convexity.

\begin{definition}\label{Def_GammaConvex}
A subset $Y$ of $\MS$ is called \emph{$\WT$-convex} if it is the intersection of finitely many special half-apartments. The \emph{$\WT$-convex hull} $c_{\WT}(Y)$ of a subset $Y$ of $X$ is the intersection of all special half-apartments containing $Y$.
\end{definition}

Note that Weyl chambers and hyperplanes are $\aW$-convex, as well as finite intersections of $\aW$-convex sets. Special hyperplanes and Weyl chambers are also $\WT$-convex.

\begin{prop}{\cite[Prop.2.13]{BennettDiss}}\label{Prop_segment}\index{segment}
For any two special vertices $x,y$ in the model space $\MS$ the \emph{segment} $\seg(x,y)=\{z\in\MS : d(x,y)=d(x,z)+d(z,y)\}$ is the same as the $\WT$-convex hull $c_{\WT}(\{x,y\})$.
\end{prop}

\begin{definition}\label{Def_parallel}
Two subsets $\Omega_1,\Omega_2$ of a $\Lambda$-metric space are at \emph{bounded distance}\index{bounded distance} if there exists $N\in\Lambda$ such that for all $x\in \Omega_i$ there exists $y\in \Omega_j$ such that $d(x,y)\leq N$ for $\{i,j\}=\{1,2\}$.
Subsets of a metric space are \emph{parallel} if they are at bounded distance. 
\end{definition}

Note that parallelism is an equivalence relation. One can prove

\begin{prop}{\cite[Section 2.4]{Bennett}}
Let $\MS=\MS(\RS, \Lambda)$ equipped with the full affine Weyl group $\aW$. Then the following is true
\begin{enumerate}
\item Two hyperplanes or Weyl simplices are parallel if and only if they are translates of each other by elements of $\aW$.
\item For any two parallel Weyl chambers $S$ and $S'$ there exists a Weyl chamber $S''$ contained in $S\cap S'$ and parallel to both. 
\end{enumerate}
\end{prop}

Moreover $H_{\alpha, k}$ is parallel to $t_{\beta^\vee}(H_{\alpha, k})=H_{\alpha, k+\beta^\vee(\alpha)},$ for all $\beta\in\RS$
where $t_{\beta^\vee}$ is the translation in $\MS$ by $\frac{2}{(\beta, \beta)}\beta$.

%
%
%
%
%

\section{The definition of generalized affine buildings}\label{Sec_generalizedAffineBuildings}

Throughout the following let $\MS=\MS(\RS,\Lambda, T)$ be as defined in \ref{Not_modelSpace}.

\begin{definition}\label{Def_LambdaBuilding}
Let $X$ be a set and $\App$ a collection of injective charts $f:\MS\hookrightarrow X$.
We call the images $f(\MS)$  of charts $f\in\App$ \emph{apartments} of $X$ and we define \emph{Weyl chambers, hyperplanes, half-apartments, special vertices, ... of $X$} to be images of such in $\MS$ under any $f\in\App$. The set $X$ is a \emph{(generalized) affine building} with \emph{atlas} (or \emph{apartment system}) $\App$ if the following conditions are satisfied
\begin{enumerate}[label={(A*)}, leftmargin=*]
\item[(A1)] For all $f\in\App$ and $w\in \WT$ the concatenation $f\circ w$ is contained in $\App$. 
\item[(A2)] Given two charts $f,g\in\App$ with $f(\MS)\cap g(\MS)\neq\emptyset$. Then $f^{-1}(g(\MS))$ is a closed convex subset of $\MS$. There exists $w\in \WT$ with $f\vert_{f^{-1}(g(\MS))} = (g\circ w )\vert_{f^{-1}(g(\MS))}$.
\item[(A3)] For any two points in $X$ there is an apartment containing both.
\end{enumerate}
Axioms $(A1)-(A3)$ imply the existence of a $\Lambda$-distance on $X$, that is a function $d:X\times X\mapsto \Lambda$ satisfying all conditions of Definition~\ref{Def_metric} but the triangle inequality.  Define the distance of points $x,y$ in $X$ to be the distance of their preimages under a chart $f$ of an apartment containing both.
\begin{enumerate}[label={(A*)}, leftmargin=*]
\item[(A4)] Given Weyl chambers $S_1$ and $S_2$ in $X$ there exist sub-Weyl chambers $S_1', S_2'$ in $X$ and $f\in\App$ such that $S_1'\cup S_2' \subset f(\MS)$. 
\item[(A5)] For any apartment $A$ and all $x\in A$ there exists a \emph{retraction} $r_{A,x}:X\to A$ such that $r_{A,x}$ does not increase distances and $r^{-1}_{A,x}(x)=\{x\}$.
\item[(A6)] Let $f, g$ and $h$ be charts such that the associated apartments pairwise intersect in half-apartments. Then $f(\MS)\cap g(\MS)\cap h(\MS)\neq \emptyset$. 
\end{enumerate}
By $(A5)$ the distance function $d$ on $X$ is well defined and satisfies the triangle inequality.

The \emph{dimension} of the building $X$ is $n=\rk(\RS)$, where $\MS\cong(\Lambda)^n$. 
\end{definition}

Tits defined his ``syst{\`e}me d'appartements'' in \cite{TitsComo} by giving five axioms. The first four are the same as $(A1)-(A4)$ above. The fifth axiom originally reads different from ours but was later replaced with $(A5)$ as presented in the definition above. In fact if $\Lambda=\R$ axiom $(A6)$ follows from $(A1)-(A5)$. But in the general case this additional axiom is necessary as illustrated with an example given on p. 563 in \cite{Bennett}.
However in \cite{Bennett}  axiom $(A6)$ is mostly used to avoid pathological cases and to guarantee the existence of the panel and wall trees.
One can find a short history of Tits' axioms in \cite{Ronan}. Equivalent sets of axioms are discussed by Parreau in \cite{Parreau}.

\begin{definition}\label{Def_thick}
\index{generalized affine building!thick}
Let $X$ be a generalized affine building with model space $\MS(\RS, \Lambda, T)$ and apartment system $\App$. We call  $X$  \emph{thick with respect to $\WT$} if for any special hyperplane $H$ of $X$ there exist apartments $A_1=f_1(\MS)$ and $A_2=f_2(\MS)$, with $f_i\in\App, i=1,2$ such that $H\in A_i, i=1,2$ and $A_1\cap A_2$ is one of the two half-apartments of $A_1$ (or $A_2$) determined by $H$.  Furthermore apartments do not branch at non-special hyperplanes.
\end{definition}

\begin{remark}
If in the previous definition $T=\MS$ then $X$ is a building branching everwhere.
\end{remark}

\begin{definition}\label{Def_iso}
\index{{generalized affine building}!{isomorphism}}
Two affine buildings $(X_1,\App_1),(X_2,\App_2)$ of the same type $\MS(\Lambda, \RS)$ are \emph{isomorphic} if there exist maps $\pi_1: X_1 \to X_2$, $\pi_2: X_2 \to X_1$, further maps $\pi_{\App_1}:\App_1 \to \App_2$ , $\pi_{\App_1}:\App_2 \to \App_1$ and an automorphism $\sigma$ of $\MS$ such that 
\begin{align*}
\pi_i\circ\pi_j=\one_{X_i} & \text{ with } \{i,j\}=\{1,2\},\\
\pi_{\App_i}\circ \pi_{\App_i}=\one_{\App_i} & \text{ with } \{i,j\}=\{1,2\},
\end{align*}
and the following diagram commutes for all $f\in\App_i$ with $\{i,j\}=\{1,2\}$
\[ \begin{xy}
 	\xymatrix{
		\MS \ar[d]_\sigma \ar[r]^f  & X_i \ar[d]^{\pi_i} \\
		\MS \ar[r]_{\pi_{\App_{i}}(f)} & X_j
	}.
\end{xy}\]
\end{definition}

Examples of a generalized affine buildings are $\Lambda$-trees without leaves. The definition of a $\Lambda$-tree, \cite[p.560]{Bennett}, is equivalent to the definition of an affine building of dimension one.
Generalized affine buildings arise for example from groups defined over fields with $\Lambda$-valued valuations. An example is given in \cite[Example 3.2]{Bennett} associating a generalized affine buildings to $SL_n(K)$ where $K$ is a field with $\Lambda$-valued valuation.

Note that the Davis realization of a simplicial affine building is a generalized affine building, as defined in~\ref{Def_LambdaBuilding} with $\Lambda=\R$ and $T$ chosen equal to the co-root lattice $\QQ(\RS^\vee)$ of $\RS$.

\section{Local and global structure}\label{Sec_local-global}

Any simplicial affine building has an associated spherical building at infinity. This useful and important result by Bruhat and Tits \cite{BruhatTits} is also true in the generalized setting.

\begin{definition}\label{Def_buildingAtInfinity}
\index{building at infinity}
Let $(X,\App)$ be an affine building. We denote by $\partial S$ the parallel class of a Weyl chamber $S$ in $X$. Let
$$\binfinity X =\{\partial S : S \text{ Weyl chamber of } X \text{ contained in an apartment of } \App\}$$
be the set of chambers of the \emph{spherical building at infinity} $\binfinity X$. We say that two chambers $\partial S_1$ and $\partial S_2$ are \emph{adjacent} if there exist representatives $S'_1, S'_2$ contained in a common apartment having  the same basepoint and are adjacent in $X$.
\end{definition}

\begin{prop}\label{Prop_buildingAtInfinity}
Let $(X,\App)$ be an affine building modeled on $\MS(\RS, \Lambda, T)$. The set $\binfinity X$ defined above is a spherical building of type $\RS$ with apartments in one to one correspondence with apartments of $X$.
\end{prop}
\begin{proof}
It is obvious that $\binfinity X$ is a simplicial complex with adjacency as defined in \ref{Def_buildingAtInfinity}. 
An apartment in $\binfinity X$ is defined to be the set of equivalence classes determined by Weyl chambers in an apartment of $X$. One easily observes that these are Coxeter complexes of type $\RS$ and that hence $\binfinity X$ has to be of this type.

Given two chambers $c$ and $d$ in $\binfinity X$, let $S$ and $T$ be representatives of $c$, respectively $d$. By axiom $(A4)$ there exists an apartment $A$ containing sub-Weyl chambers of $S$ and $T$. The set of equivalence classes of Weyl chambers determined by $A$ hence contains $c$ and $d$. Therefore any two chambers are contained in a common apartment, hence 2. of the definition of a building as given on p. 77 in \cite{Brown} holds.

If $\partial A$ and $\partial A'$ are two apartments of $\binfinity X$ both containing the chambers $c$ and $d$, then there exist charts $f$ and $f'$ such that $f(\MS)$ and $f'(\MS)$ contain representatives $S, S'$ of $c$ and $T, T'$ of $d$. The Weyl chambers $S,S'$ and $T,T'$ intersect in sub-Weyl chambers $S''$ and $T''$, respectively. The map $f'\circ f^{-1}$ fixes $S''$ and $T''$ and induces an isomorphism from $\partial A$ to $\partial A'$. Therefore $\binfinity X$ is indeed a spherical building.
\end{proof}

In contrast to a remark made in \cite{BennettDiss} it is possible to prove Proposition \ref{Prop_buildingAtInfinity} without using axiom $(A5)$.

The local structure of an affine building was not examined in \cite{BennettDiss}. In analogy to the residues of vertices in a simplicial affine building one can associate to a vertex of a generalized affine building a spherical building. Most of the following in based on \cite{Parreau}.

Let in the following $(X,\App)$ be an affine building of type $\MS=\MS(\Lambda, \RS, T)$ and let $\binfinity X$ denote its spherical building at infinity.

\begin{definition}\label{Def_germ}
\index{{generalized affine building}!{germ}}
Two Weyl simplices $S$ and $S'$ \emph{share the same germ} if both are based at the same vertex and if $S\cap S'$ is a neighborhood of $x$ in $S$ and in $S'$.
It is easy to see that this is an equivalence relation on the set of Weyl simplices based at a given vertex. The equivalence class of $S$, based at $x$, is denoted by $\Delta_x S$ and is called the \emph{germ of $S$ at $x$}.
\end{definition}

\begin{remark}\label{Rem_partialorder}
The germs of Weyl simplices at a special vertex $x$ are partially ordered by inclusion: $\Delta_x S_1$ is contained in $\Delta_xS_2$ if there exist $x$-based representatives $S'_1, S'_2$ contained in a common apartment such that $S_1'$ is a face of $S_2'$. Let $\Delta_xX$ be the set of all germs of Weyl simplices based at $x$.
\end{remark}

We say that a germ $\mu$ of a Weyl chamber $S$ at $x$ is contained in a set $Y$ if there exists $\varepsilon\in\Lambda^{+}$ such that $S\cap B_\varepsilon(x)$ is contained in $Y$.

\begin{prop}\label{Prop_tec16}
Let $(X, \App)$ be an affine building and $c$ a chamber in $\binfinity X$. Let $S$ be a Weyl chamber in $X$ based at $x$. Then there exists an apartment $A$ such that $\Delta_xS \subset A$ and  $c\in \partial A$.
\end{prop}
The proof of the above proposition is precicely the same as the proof of Proposition 1.8 in \cite{Parreau}.

\begin{corollary}\label{Cor_WeylChamber}
Fix a point $x\in X$. For each Weyl simplex $F$ there exists a unique Weyl simplex $F'$ of the same dimension based at $x$ and parallel to $F$.
\end{corollary}
\begin{proof}
Apply Proposition~\ref{Prop_tec16} to $x$ and $c=\partial F$ and arbitrary $S$ based at $x$.
\end{proof}

\begin{corollary}\label{Cor_tec17}
For any chamber $c\in\binfinity X$ the affine building $X$ is as a set the union of all apartments containing a representative of $c$.
\end{corollary}
\begin{proof}
Fix a chamber $c$ at infinity. For all  points $x\in X$ and arbitrary Weyl chambers $S$ based at $x$ there exists by \ref{Prop_tec16} an apartment $A$ such that $A$ contains $x$ and a germ of $S$ at $x$ and such that $c$ is contained in $\partial A$. 
\end{proof}

\begin{corollary}\label{Cor_tec18}
Given a germ $\mu$ of a Weyl chamber. Then $X$ is the union of all apartments containing $\mu$.
\end{corollary}
\begin{proof}
Assume $\mu=\Delta_xS$. For all $y\in X$ there exists, by axiom $(A3)$, an apartment $A'$ containing $x$ and $y$. Let $T$ be a Weyl chamber in $A'$ based at $x$ containing $y$ and let $c=\partial T$. By Proposition~\ref{Prop_tec16} there exists an apartment $A$ such that a germ of $S$ at $x$ is contained in $A$ and such that the corresponding apartment $\partial A$ of $\binfinity X$ contains $c$. But then the unique representative $T$ of $c$ based at $x$ is also contained in $A$. Therefore $A$ contains $y$ and a germ of $S$ at $x$.
\end{proof}

\begin{corollary}\label{Cor_GG}
Any two germs of Weyl chambers based at the same vertex are contained in a common apartment, that is two chambers of $\Delta_xX$ are contained in a common apartment of $\Delta_xX$.
\end{corollary}
\begin{proof}
Let $S$ and $T$ be Weyl chambers both based at $x$. By Proposition~\ref{Prop_tec16} there exists an apartment $A$ of $X$ containing $S$ and a germ of $T$ at $x$. Therefore $\Delta_xS$ and $\Delta_xT$ are both contained in the apartment $\Delta_xA$.
\end{proof}

\begin{prop}\label{Prop_A3'}
Let $(X,\App)$ be an affine building. Let $S$ and $T$ be Weyl chambers based at $x$ and $y$, respectively. Then there exists an apartment $A$ of $X$ containing a germ of $S$ at $x$ and a germ of $T$ at $y$.
\end{prop}
\begin{proof}
By axiom $(A3)$ there exists an apartment $A$ containing $x$ and $y$. We choose an $x$-based Weyl chamber $S_{xy}$ in $A$ that contains $y$ and denote by $S_{yx}$ the Weyl chamber based at $y$ such that $\partial S_{xy}$ and $\partial S_{yx}$ are opposite in $\partial A$. Then $x$ is contained in $S_{yx}$. If $\Delta_yT$ is not contained in $A$ apply Proposition~\ref{Prop_tec16} to obtain an apartment $A'$ containing a germ of $T$ at $y$ and containing  $\partial S_{yx}$ at infinity. But then $x$ is also contained in $A'$. 

Let us denote by $S'_{xy}$ the unique Weyl chamber contained in $A'$ having the same germ as $S_{xy}$ at $x$.
Without loss of generality we may assume that the germ $\Delta_yT$ is contained in $S'_{xy}$. Otherwise $y$ is contained in a face of $S'_{xy}$ and we can replace $S'_{xy}$ by an adjacent Weyl chamber in $A'$ satisfying this condition.
A second application of Proposition~\ref{Prop_tec16} to $\partial S'_{xy}$ and the germ of $S$ at $x$ yields an apartment $A''$ containing $\Delta_xS$ and $S'_{xy}$ and therefore $\Delta_yT$.
\end{proof}

\begin{remark}
Corollaries \ref{Cor_WeylChamber} and \ref{Cor_GG} are the direct analogs of 1.9 and 1.11 of \cite{Parreau} and Proposition \ref{Prop_A3'} corresponds to \cite[1.16]{Parreau}. Note, however, that the proof is different.
\end{remark}

\begin{defthm}\label{Thm_residue}
\index{{generalized affine building}!{residue}}
Let $(X,\App)$ be an affine building with model space $\MS(\RS,\Lambda,T)$. Then $\Delta_xX$ is a spherical building of type $\RS$ for all $x$ in $X$. If $x$ is special and $X$ is thick with respect to $\WT$, then $\Delta_xX$ is thick as well. Furthermore $\Delta_xX$ is independent of $\App$.
\end{defthm}
\begin{proof}
We verify the axioms of the definition of a simplicial building, which can be found on  page 76 in \cite{Brown}.
It is easy to see that $\Delta_xX$ is a simplicial complex with the partial order defined in \ref{Rem_partialorder}. It is a pure simplicial complex, since each germ of a face is contained in a germ of a Weyl chamber. The set of equivalence classes determined by a given apartment of $X$ containing $x$ is a subcomplex of $\Delta_xX$ which is, obviously, a Coxeter complex of type $\RS$. Hence we define those to be the apartments of $\Delta_xX$. Therefore, by definition, each apartment is a Coxeter complex. 
Two apartments of $\Delta_xX$ are isomorphic via an isomorphism fixing the intersection of the corresponding apartments of $X$, hence fixing the intersection of the apartments of $\Delta_xX$ as well. Finally due to Corollary~\ref{Cor_GG} any two chambers are contained in a common apartment and we can conclude that $\Delta_xX$ is a spherical building of type $\RS$.

Assume that $x$ is special and $X$ is thick with respect to $\WT$. Let $c$ be a chamber in $\Delta_xX$ and $\Delta_xA$ an apartment containing $c$. For each panel $p$ of $c$ there exists a chamber $c'$ contained in $\Delta_xA$ such that $c\cap c'=p$. The panel $p=\Delta_xF$ determines a wall $H\subset A$. Since $X$ is thick there exists an apartment $A'$ whose intersection with $A$ is a half-apartment bounded by $H$. Hence there is a third chamber $c''$ of $\Delta_xX$ determined by a Weyl chamber in $A'$ based at $x$ containing $F$. Therefore $\Delta_xX$ is thick.

Let $\App'$ be a different system of apartments of $X$ and assume w.l.o.g.~that $\App\subset \App'$.
We will denote by $\Delta$ the spherical building of germs at $x$ with respect to $\App$ and  by $\Delta'$ the building at $x$ with respect to $\App'$. Since spherical buildings have a unique apartment system $\Delta$ and $\Delta'$ are equal if they contain the same chambers. Assume there exists a chamber $c\in\Delta'$ which is not contained in $\Delta$. Let $d$ be a chamber opposite $c$ in $\Delta'$ and $a'$ the unique apartment containing both. Note that $a'$ corresponds to an apartment $A'$ of $X$ having a chart in $\App'$. There exist $\App'$-Weyl chambers $S_c$, $S_d$ contained in $A$ representing $c$ and $d$, respectively. Choose a point $y$ in the interior of $S_c$ and let $z$ be contained in the interior of $S_d$. By axiom $(A3)$ there exists a chart $f\in \App$ such that $A=f(\MS)$ contains $y$ and $z$. The apartment $A$ also contains $x$ since $x\in\seg(y,z)$ and $\seg(y,z)\subset A\cap A'$. By construction the unique Weyl chambers of $A$ based at $x$ containing $y$, respectively $z$, have germs $c$, respectively $d$, which is a contradiction. Hence $\Delta=\Delta'$.
\end{proof}

\begin{remark}
Let $\MS=\MS(\RS,\Lambda,T)$ be the model space of an affine building and let $\partial \MS$ be canonically identified with the associated spherical Coxeter complex. Note that for each $x\in \MS$ and each chamber $c$ of $\partial \MS$ there exists a Weyl chamber $S$ contained in $c$ and based at $x$. Therefore the type of $\Delta_xX$ for an affine building $(X,\App)$ modeled on $\MS$ is always $\RS$. 

If $(X,\App)$  is the geometric realization of a simplicial affine building then $\Delta_xX$ is canonically isomorphic to the residue (or link) of $x$ if and only if $x$ is a special vertex. The definition of a spherical building corresponding to the residue of a non-special vertex would be possible defining a second class of Weyl chambers based at a vertex $x$ with respect to the stabilizer $(\WT)_x$ of $x$ in the restricted affine Weyl group $\WT$. Since we will not make use of this fact, we will not give details here.
\end{remark}

\begin{prop}~\label{Prop_epi}
Let $(X,\App)$ be an affine building, $\binfinity X$ its building at infinity. For all vertices $x\in X$ there exists an epimorphism 
$$\pi_x:\; \partial_\App X \rightarrow \Delta_xX .$$
\end{prop}
\begin{proof}
Given $x\in X$ and $c\in \binfinity X$. Let $S$ be the Weyl chamber based at $x$ and contained in $c$, which exists by Corollary~\ref{Cor_WeylChamber}. Define $\pi_x(c)=\Delta_xS$, the germ of $S$ at $x$. Since for all $d \in \Delta_xX$ there exists a Weyl chamber $S'$ in $X$ such that $\Delta_xS' = d$ the map $\pi_x$ is surjective. By definition of $\Delta_xX$, the partial order and adjacency is preserved.
\end{proof}

\begin{prop}\label{Prop_A5}
Let $X$ be an affine building and let $A_i$ with $i=1,2,3$ be three apartments of $X$ pairwise intersecting in half-apartments. Then $A_1\cap A_2\cap A_3$ is either a half-apartment or a hyperplane.
\end{prop}
\begin{proof}
For $i\neq j$ denote the intersection $A_i\cap A_j$ by $M_{ij}$. The corresponding apartments $\partial A_i$ in the spherical building at infinity do as well intersect in half-apartments. Hence $\partial A_1\cap \partial A_2 \cap \partial A_3$ is either a half-apartment or a hyperplane in $\binfinity X$. Assume that $\partial A_1\cap \partial A_2 \cap \partial A_3$ is a half-apartment. Then $A_1\cap A_2\cap A_3$ is a half-apartment contained in each of the $A_i$.

Assume now that we are in the case where $\partial A_1\cap \partial A_2 \cap \partial A_3$ is a hyperplane $m$ in $\binfinity X$.
Walls at infinity correspond to parallel classes of hyperplanes in the affine building. Hence there are three hyperplanes $H_{ij}$ bounding the half-apartments $M_{ij}= A_i\cap A_j$ which are all contained in $m$. Note that the half-apartments $M_{13}$ and $M_{23}$ are opposite in $A_3$ in the sense that their union equals $A_3$. By axiom $(A6)$ the intersection $A_1\cap A_2\cap A_3$ is nonempty and equal to the strip $M_{13}\cap M_{23}$. It is obvious that the hyperplanes $H_{13}$ and $H_{23}$ are contained in $M_{13}\cap M_{23}$. This argument is symmetric in the indices. Thus $H_{12}$ is contained in $M_{13}\cap M_{23}$. Again by symmetry each of the hyperplanes is between the other two and hence $H_{12}=H_{13}=H_{23}= A_1\cap A_2 \cap A_3$.
\end{proof}

This leads to the following observation.

\begin{figure}[htbp]
 \begin{center}
 	\resizebox{!}{0.2\textheight}{\input{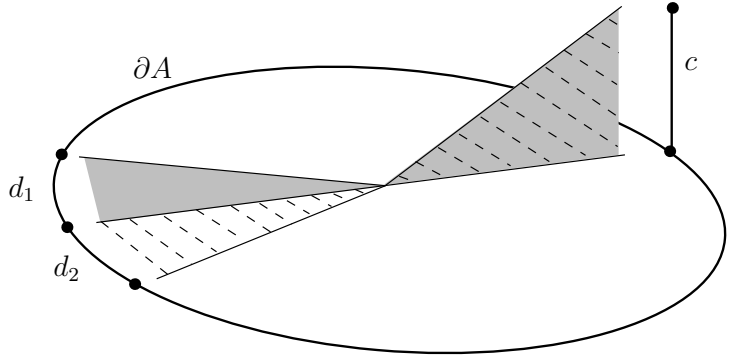}}
 	\caption[adjacent]{The sundial configuration.}
 	\label{Fig_sundial}
 \end{center} 
 \end{figure}

\begin{property}{\bf The sundial configuration. } \label{Prop_sundial}
Let $A$ be an apartment in $X$ and let $c$ be a chamber not contained in $\partial A$ but containing a panel of $\partial A$. Then $c$ is opposite to two uniquely determined chambers $d_1$ and $d_2$ in $\partial A$. Hence there exist apartments $A_1$ and $A_2$ of $X$ such that $\partial A_i$ contains $d_i$ and $c$ with $i=1,2$.  The three apartments $\partial A_1,\partial A_2$ and $\partial A$ pairwise intersect in half-apartments. Axiom $(A6)$ together with the proposition above implies that their intersection is a hyperplane. Compare Figure~\ref{Fig_sundial}.
\end{property}

\begin{prop}\label{Prop_liftGallery}
Let $x$ be an element of $X$. Let $(c_0, \ldots, c_k)$ be a minimal gallery in $\binfinity X$. We denote by $S_i$ the $x$-based representative  of $c_i$. If $(\pi_x(c_0), \ldots, \pi_x(c_k))$ is minimal in $\Delta_xX$, then there exists an apartment containing $\cup_{i=0}^k S_i$.
\end{prop}
\begin{proof}
The proof is by induction on $k$. For $k=0$ there is just one Weyl chamber and the result holds. Let $A'$ be an apartment containing $S_1\cup S_2\cup \ldots \cup S_{k-1}$. If $c_k$ is contained in $\partial A'$ we are done. If $c_k$ is not contained in $\partial A'$ we have the sundial configuration, which determines a unique hyperplane $H$ in $A'$. Let $H^+$ be the unique half-apartment of $A'$ determined by this hyperplane which contains a representative of $c_0$. Then $H^+$  also contains representatives of $c_1,\ldots,c_{k-1}$, since this is a minimal gallery and $S_k$ is on the other side of $H$. 

We claim that $x$ is contained in $H^+$.
Let $A''$ be the apartment in the sundial configuration containing $H^- = (A'\setminus H^+)\cup H$ and $c_k$ at infinity. If $x$ is contained in $A'\setminus H^+ \subset A''$ the Weyl chamber $S_k$ is contained in $A''$. Let $\rho: A'' \rightarrow A'$ be the isometry  fixing $A''\setminus H^+$. The Weyl chamber $S_k$ is mapped onto $S_{k-1}$ and the set $S_k\cap(A'\cap H^+)$ is pointwise fixed. Therefore $\pi_x(c_{k-1})=\pi_x(c_k)$ which is a contradiction. Hence $x$ is contained in $H^+$ and $\cup_{i=0}^{k-1} S_i \subset H^+$. Let $A$ now be the apartment in the sundial configuration containing $H^+$ and $S_k$. Then $\cup_{i=0}^{k} S_i$ is contained in $A$.
\end{proof}

\begin{remark}
Propositions \ref{Prop_A5} and \ref{Prop_liftGallery} as well as \ref{Prop_sundial} are due to Linus Kramer.
\end{remark}

\begin{corollary}\label{Cor_CO}
If the germs of two Weyl chambers based at the same vertex are opposite in the corresponding residue, then these Weyl chambers are contained in a unique common apartment. 
\end{corollary}
\begin{proof}
Choose a minimal gallery $(c_0,c_1,\ldots, c_n)$ from $c_0=\partial S$ to $c_n=\partial T$ and consider the representatives $S_i$ of $c_i$ based at $x$. Then $S_0=S$ and $S_n=T$. Proposition~\ref{Prop_liftGallery} implies the assertion.
\end{proof}

\begin{thm}\label{Thm_projection1}
Let $(X,\App)$ be an affine building, $A$ an apartment and $c,d$ opposite chambers in $\partial A$. Then
$$
x\in A \Longleftrightarrow \pi_x(c) \text{ and } \pi_x(d) \text{ are opposite in } \Delta_xX. 
$$ 
The restriction of $\pi_x$ to the boundary of an apartment $A$ containing $x$ is an isomorphism onto its image.
\end{thm}
\begin{proof}
First we assume that $x\in A$. Each panel $p_i$ of $c$ defines an equivalence class $m_i$ of parallel hyperplanes in $A$. The classes $m_i$ are hyperplanes in $\partial A$. Let $H_i$ be the unique representative of  $m_i$ in $A$ containing $x$. Let further $M_i$ be the (unique) half-apartment of $A$ determined by $H_i$ such that $\partial M_i$ contains $c$. The intersection of all $M_i$ is the Weyl chamber $S$ based at $x$, representing $c$. Similarly we have a Weyl chamber $T$ based at $x$ representing $d$. The defining walls of $S$ and $T$ are the same. Thus the chambers $\pi_x(c)$ and $\pi_x(d)$ are opposite in $\Delta_xX$.

Given $c,d$ in $\binfinity X$ such that the chambers $\pi_x(c)$ and $\pi_x(d)$ are opposite. Let $S$ and $T$ be the Weyl chambers based at $x$ contained in $c$, $d$, respectively.
Choose a minimal gallery $\gamma'=(c'_0=\pi_x(c), c'_1, \ldots, c'_{n-1}, c'_n=\pi_x(d))$ in the residue $\Delta_xX$ at $x$.  Then there exists a minimal gallery $\gamma=(c_0=c, c_1, \ldots, c_{n-1}, c_n=d)$ in $\binfinity X$ such that $\pi_x(c_i)=c'_i$. Denote by $S_i$ the unique Weyl chamber based at $x$ and contained in $c_i$. Proposition~\ref{Prop_liftGallery} implies the existence of an apartment $A$ containing $\cup_{i=0}^k C_i$ and hence $x$. Uniqueness is clear by $(A2)$. 
\end{proof}

%
%
%
%
%

\section{Retractions}\label{Sec_retractions}

We will define two types of retractions. On the one hand a retraction centered at a germ of a Weyl chamber and on the other hand a retraction centered at a chamber in the spherical building at infinity.

\begin{definition}\label{Def_vertexRetraction}
\index{vertex retraction}
Let $(X,\App)$ be an affine building. We fix a Weyl chamber $S$ based at a vertex $x$ in $X$ and we denote its germ by $\Delta_xS$. According to Corollary~\ref{Cor_tec18} the building $X$ is (as a set) the union of all apartments containing $\Delta_xS$.
Given $y\in X$ choose a chart $g\in\App$ such that $y$ and $\Delta_xS$ are contained in  $g(\MS)$ and define 
$$ r_{A,\Delta_xS}(y) = (f\circ w\circ g^{-1} )(y)  $$
where $w\in\aW$ is such that $g\vert_{g^{-1}(f(\MS))}=(f\circ w)\vert_{g^{-1}(f(\MS))}$.
The map $r_{A,\Delta_xS}$ is called \emph{retraction onto $A$ centered at $\Delta_xS$}.
\end{definition}

\begin{definition}\label{Def_retractionInfty}
\index{retraction centered at infinity}
Let $(X,\App)$ be an affine building. we fix an equivalence class of Weyl chambers $c \subset \binfinity X$ and an apartment $A=f(\MS)$ containing some representative $S$ of $c$. By Corollary~\ref{Cor_tec17} the building $X$ is the union of apartments containing a sub-Weyl chamber of $S$. For $x\in X$ we choose a chart $g\in\App$ such that $x\in g(\MS)$ and such that $c\in\partial (g(\MS))$. We define 
$$
 \rho_{A,c}(x) = (f\circ w\circ g^{-1} )(x)  
$$
where $w\in\aW$ is such that $g\vert_{g^{-1}(f(\MS))}=(f\circ w)\vert_{g^{-1}(f(\MS))}$.
The map $\rho_{A,c}$ is called \emph{retraction onto $A$ centered at $c=\partial S$} or (being slightly imprecise) \emph{retraction centered at infinity}.
\end{definition}

\begin{prop}\label{Prop_r-rho}
Let $(X,\App)$ be an affine building. Fix an apartment $A$ of $X$. Let $S$ be a Weyl chamber contained in $A$ and let $c$ be a chamber in $\partial A$. Then the following hold:
\begin{enumerate}
  \item The maps $r_{A,\Delta_xS}$ and $\rho_{A,c}$ are well defined.
  \item \label{tec123} The restriction of $\rho_{A,c}$ to an apartment $A'$ containing $c$ at infinity is an isomorphism onto $A$.
  \item The restriction of the retraction $r_{A,\Delta_xS}$ to an apartment $A'$ containing $\Delta_xS$ is an isomorphism onto $A$.
\end{enumerate}
\end{prop}
\begin{proof}
The second and third assertions are clear by definition.

Assume $A_i\define f_i(\MS)$, $i=1,2$, are two apartments both containing $\Delta_xS$ and a point $y$. We let $w_i$ be the element of $\aW$ appearing in the definition of $r_{A,\Delta_xS}(y)$ with respect to $f_i$. It suffices to prove
\begin{equation}\label{Equ_tec28}
f\circ w_1\circ f_1^{-1}(y)=f\circ w_2\circ f_2^{-1}(y).
\end{equation}
By assumption the germ $\Delta_xS$ is contained in $A_1\cap A_2$ hence there exists by $(A2)$ an element $w_{12}\in\aW$ such that 
$$ 
f_2\circ w_{12}\;\vert_{f_1^{-1}(f_2(\MS))} = f_1 \;\vert_{f_1^{-1}(f_2(\MS))} .
$$
Since $y\in A_1\cap A_2$, we have 
\begin{equation}\label{Equ_tec29}
f\circ w_2\circ f_2^{-1} 
	= f\circ w_2\circ f_2^{-1}(f_2\circ w_{12}(f_1^{-1}(y)))
	= f\circ w_2 w_{12}(f_1^{-1}(y)).
\end{equation}
There are unique Weyl chambers $S_1$ and $S_2$ contained in $A_1$ and $A_2$, respectively, satisfying the property that $\Delta_xS_i=\Delta_xS$, $i=1,2$. Since equation (\ref{Equ_tec29}) is true for all $y\in A_1\cap A_2$, it is in particular true for the intersection $C$ of the Weyl chambers $S_1$ and $S_2$. Therefore 
$$
f\circ w_1\circ f_1^{-1}(C) = f\circ w_2\circ w_{12} \circ f_1^{-1}(C)
$$
and hence $w_2 w_{12}=w_1.$ Combining this with (\ref{Equ_tec29}) yields equation (\ref{Equ_tec28}).

An argument along the same lines proves that $\rho_{A,c}$ is well defined.
\end{proof}

\begin{lemma}\label{Lem_compRetractions}
Let $(X,A)$ be an affine building.
\begin{enumerate}
\item Given a Weyl chamber $S$ in $X$ and apartments $A_i$, where $i$ ranges from $1$ to $n$, containing sub-Weyl chambers of $S$. Denote by $\rho_i$ the retraction $\rho_{A_i, \partial S}$. Then
$$(\rho_1\circ\rho_2\circ\ldots\circ\rho_n) = \rho_1 .$$
\item Let $\Delta_xS$ be a germ of a Weyl chamber $S$ at $x$. Let $A_i, i=1,\ldots, n$ be a set of apartments containing $\Delta_xS$ and denote by $r_i$ the retraction onto $A_i$ centered at $\Delta_xS$. Then 
$$(r_1\circ r_2\circ\ldots\circ r_n) = r_1 .$$
\end{enumerate}
\end{lemma}
\begin{proof}
According to \ref{Prop_r-rho} (\ref{tec123}) the restriction of $\rho_i$ to an apartment containing a sub-Weyl chamber of $S$ is an isomorphism for all $i$. By Corollary~\ref{Cor_tec17} the building $X$ is as a set the union of all apartments containing a representative of $\partial S$. Therefore (1) follows. Similar arguments using Corollary~\ref{Cor_tec18} imply the second part of the lemma.
\end{proof}

\begin{prop}\label{Prop_compRetractions}
Let $(X,\App)$ be an affine building modeled on $\MS(\RS,\Lambda)$. For all retractions $\rho_{A,c}$ centered at infinity and all $x\in X$ there exists $y\in A$ such that $$\rho_{A,c}(x)=r_{A,\pi_y(c)}(x),$$ where $\pi_y$ is as in Proposition~\ref{Prop_epi}.
\end{prop}
\begin{proof}
By Proposition~\ref{Prop_tec16} there exists an apartment $B$ of $X$ containing $x$ and a Weyl chamber $S$ contained in $c$. The intersection of $A$ and $B$ contains a sub-Weyl chamber $S'$ of $S$ and the restriction of $\rho_{A,c}$ to $B$ is an isomorphism onto $A$ fixing $A\cap B$ pointwise. Let $y$ be a vertex in $S'$ and denote by $S_y$ the Weyl chamber based at $y$ and parallel to $S'$. The restriction of $r_{A,\Delta_yS}$ to $B$ is by \ref{Prop_r-rho} an isomorphism onto $A$ which fixes $A\cap B$ pointwise. Therefore $\rho_{A,c}(x)=r_{A,\Delta_yS}(x)$.
\end{proof}

%
%
%
%
%
%

\section{Finite covering properties}\label{Sec_FC}

In the present section we will introduce certain finite covering properties which are (in Section \ref{Sec_convexityRevisited}) used to prove that the retractions defined at the end of the previous section are distance non-increasing.

We first prove that segments are contained in apartments. As always let $(X,\App)$ be a generalized affine building. Recall the definition of the \emph{segment} of points $x$ and $y$ in $X$:
$$\seg(x,y)\define\{z\in X : d(x,y)=d(x,z)+d(z,y) \}.$$

\begin{lemma}\label{Lem_tec31}
Let $x,y$ be points in $X$ and let $A$ be an apartment containing $x$ and $y$. Let $r$ be a retraction onto $A$ existing by axiom $(A5)$. Then for all $z\in \seg(x,y)$ the following is true
\begin{enumerate}
  \item $r(z)\in \seg(x,y)$ and 
  \item $d(x,z)=d(x,r(z))$ as well as $d(y,z)=d(y,r(z))$.
\end{enumerate}
\end{lemma}
\begin{proof}
Since $d$ satisfies the triangle inequality and $r$ is distance non-increasing we have
$$
d(x,y)\leq d(x,r(z)) + d(y,r(z)) \leq d(x,z)+d(y,z)=d(x,y)
$$
and $r(z)$ is contained in $\seg(x,y)$.

Therefore 
\begin{equation}\label{Equ_1}
d(x,r(z))+d(r(z),y) = d(x,z)+d(z,y).
\end{equation}
The assumption that $d(x,r(z))$ is strictly smaller than $d(x,z)$ contradicts equation \ref{Equ_1}, hence $d(x,z)=d(x,r(z))$ and for symmetric reasons $d(y,z)=d(y,r(z))$.
\end{proof}

\begin{prop}\label{Prop_SegInApp}
Segments are contained in apartments.
\end{prop}
\begin{proof}
Let $x$ and $y$ be points in $X$ and let $A$ be an apartment containing $x$ and $y$. 
Assume there exists $z$ in $\seg(x,y)\!\setminus\! A$.

Let $S$ be a Weyl chamber based at $x$ containing $z$ and choose a point $p$ in $\seg(x,z)\cap S\cap A$ such that the germ of the Weyl chamber $S'$ parallel to $S$ and based at $p$ is not contained in the apartment $A$. 

On the one hand $p\in\seg(x,z)\subset\seg(x,y)$, therefore
$$
d(x,y)-d(x,p)=d(x,p)-d(p,y)-d(x,p) =d(p,y)
$$
and on the other hand, since in addition $z\in\seg(x,y)$, we obtain 
$$
d(x,y)-d(x,p)=d(p,z)+d(z,y).
$$
Thus $z$ is contained in the segment of $p$ and $y$. 

Let $S^-$ be a Weyl chamber in $A$ based at $p$, opposite $S'$ at $p$, and chosen such that it contains $x$. The image of $z$ under the retraction $r^-\define r_{A,S^-}$ onto $A$ will be denoted by $z^-$. By \ref{Cor_CO} there exists a unique apartment $A'$ containing $S'\cup S^-$, in particular $\{x,p,z\}\subset A$. The restriction of $r^-$ to $A'$ is an isomorphism onto $A$. Thus 
$$
d(p,z)=d(p,z^-).
$$

Let $r^+\define r_{A,z^-}$ be a retraction associated to $A$ and $z^-$ which exists by axiom $(A5)$. Since $(r^+)^{-1}(z^-) = \{z^-\}$ the image $z^+$ of $z$ under $r^+$ is different from $z^-$. According to \ref{Lem_tec31} and the fact that $z\in\seg(x,y)$, the point $z^+$ is contained in $\seg(x,y)$. Hence $d(x,y)=d(x,z^+) + d(z^+,y)$. 
Since $z$ is also contained in $\seg(p,y)$, we obtain from \ref{Lem_tec31}.(2) that 
$$
d(p,z^+)=d(p,z).
$$
Using the equations above, the triangle inequality and the positivity of distance functions we can calculate
\begin{align*}
d(p,z) &= d(p,z^-)\geq d(p,z^+)-d(z^-,z^+) \geq d(p,z^+)=d(p,z).
\end{align*}
This implies $d(z^-, z^+)=0$ which contradicts $z^+\neq z^-$. Therefore $\seg(x,y) \subset A$.
\end{proof}

\begin{lemma}\label{Lem_cover}
Given an apartment $A$ and a point $z$ in $X$. Then $A$ is contained in the (finite) union of all Weyl chambers based at $z$ with equivalence class contained in $\partial A$.
\end{lemma}
\begin{proof}
In case $z$ is contained in $A$ this is obvious. Hence we assume that $z$ is not contained in $A$. For all $p\in A$ there exists, by definition, an apartment $A'$ containing $z$ and $p$. Let $S_+$ be a $p$-based Weyl chamber containing $z$.  We denote by $\sigma_+$ its germ at $p$. There exists a $p$-based Weyl chamber $S_-$ in $A$ such that its germ $\sigma_-$ is opposite $\sigma_+$ in the residue $\Delta_pX$. By Corollary~\ref{Cor_CO} the Weyl chambers $S_-$ and $S_+$ are contained in a common apartment $A''$. Let $T$ be the unique representative of $\partial S_-$ in $A$ based at $z$. Since $z\in S_+$ and $\sigma_+$ and $\sigma_-$ are opposite in $\Delta_pX$, the point $p$ is contained in $T$. Therefore $T$ contains $p$ and $\partial T$ is contained in $\partial A$. To finish the proof we observe that there are only finitely many chambers in $\partial A$, and that by Corollary~\ref{Cor_WeylChamber} for each of them exists a unique representing Weyl chamber based at $z$.
\end{proof}

\begin{prop}\label{Prop_FC'}
Let $(X,\App)$ be a generalized affine building, let $x$ and $y$ be points in $X$. For all $z\in X$ the following is true:
\begin{itemize}[label={(FC')}, leftmargin=*]
  \item[$\mathrm{(FC')}$] The segment $\seg(x,y)$ of $x$ and $y$ is contained in a finite union of Weyl chambers based at $z$.
\end{itemize}
Furthermore, is $\mu$ a germ of a Weyl chamber based at $z$, then $\seg(x,y)$ is contained in a finite union of apartments containing $\mu$.
\end{prop}
\begin{proof}
By Proposition \ref{Prop_SegInApp} there is an apartment $A$ containing the segment of $x$ and $y$.  But according to Lemma \ref{Lem_cover} the apartment $A$ is contained in a finite union Weyl chambers based at $z$, hence so is $\seg(x,y)$. The remainder of the proposition follows from \ref{Prop_tec16}.
\end{proof}

\begin{prop}\label{Prop_retraction}
For all apartments $A$ and germs $\mu$ of Weyl chambers contained in $A$ the retraction $r_{A,\mu}$, as defined in \ref{Def_vertexRetraction}, is distance non-increasing.
\end{prop}
\begin{proof}
Given two points $x$ and $y$ in $X$. By \ref{Prop_FC'} there exists a finite collection of apartments $A_0,\ldots, A_n$ each containing $\mu$ such that the union contains the segment of $x$ and $y$. Let them be enumerated such that each $A_i$ contains the germ of $S$ and 
such that $A_i\cap A_{i+1} \neq \emptyset$ for all $i=0,\ldots, n-1$. 
Observe that one can find a finite sequence of points $x_i$, $i=0,\ldots,n$ with $x_0=0$ and $x_n=y$ such that 
$$
d(x,y)=\sum_{i=0}^{n-1} d(x_i,x_{i+1})
$$
and such that $A_i$ contains $x_i$ and $x_{i+1}$.
Note further that for all $i$ the restriction of $r_{A,S}$ to $A_i$ is an isomorphism onto $A$. Hence the distance $d(x_i, x_{i+1})$ of $x_i$ and $x_{i+1}$ is equal to $d(\rho(x_i), \rho(x_{i+1}))$ for all $i\neq N$. Since the metric $d$ satisfies the triangle inequality we have that $d(r(x),r(y))\leq d(x,y)$.
\end{proof}

Reading the proof of \ref{Prop_retraction} carefully it is easy to see, that besides axioms $(A1)-(A4)$ and $(A6)$ we only used the fact that the distance function on $X$, induced by the distance function on the model space, satisfies the triangle inequality. From this we can define a distance non-increasing retraction satisfying axiom $(A5)$. Therefore

\begin{corollary}
Let $(X,\App)$ be a space satisfying all axioms in Definition~\ref{Def_LambdaBuilding} but $(A5)$. Then the following are equivalent:
\begin{itemize}[label={(A5*)}, leftmargin=*]
\item[(A5)] For any apartment $A$ and all $x\in A$ there exists a \emph{retraction} $r_{A,x}:X\to A$ such that $r_{A,x}$ does not increase distances and such that $r^{-1}_{A,x}(x)=\{x\}$.
\item[(A5')] The distance function $d$ on $X$ induced by the distance function on the model space satisfies the triangle inequality.
\end{itemize}
\end{corollary}

\begin{lemma}\label{Lem_oppSectors}
Let $A$ be an apartment and let $S$ and $T$ be Weyl chambers in $A$ facing in opposite directions, that is $\partial S$ and $\partial T$ are opposite in $\partial_\App X$. If $T\cap S$ is nonempty then 
$T\cap S$ is contained in the segment of the base points $p$ of $S$ and $q$ of $T$.

Furthermore $T\cap S$ is contained in the segment of pairs of points $x\in S\!\setminus\! T$ and $y\in T\!\setminus\! S$ such that $T$ is contained in the representative $x$-based representative $T_x$ of $\partial T$ and, symmetrically, such that $S$ is contained in the $y$-based representative $S_y$ of $\partial S$. 
\end{lemma}
\begin{proof}
Since $S$ and $T$ are opposite, the intersection is the $\aW$-convex hull of the base points $p$ of $S$ and $q$ of $T$ which is by Proposition \ref{Prop_segment} equal to the segment of $p$ and $q$. Choose points $x\in S$ such that the unique Weyl chamber based at $x$ parallel to $T$ contains $T$. Analogously choose $y$. By comparison of indices of the bounding hyperplanes of $S_y$ and $T_x$ we observe that $\seg(x,y)$ contains $T\cap S$. 
\end{proof}

\begin{lemma}\label{Lem_tec40}
Let $x$, $y$ be points of $X$ and let $S$ be a Weyl chamber in an apartment $A$. 
Then there exists a sub-Weyl chamber $S'$ of $S$ such that 
$$
d(r(p),q) = d(p,q) \text{ for all } q\in S' \text{ and all } p\in \seg(x,y)
$$
where $r\define r_{A,S'}$. 
Is $z$ the basepoint of $S'$, then for all $p\in\seg(x,y)$ there exists a Weyl chamber opposite $S'$ at $z$ containing $p$.
\end{lemma}
\begin{proof}
Let $x'$ be some point in $A$ and let $\lambda$ be such that $\seg(x,y)$ is contained in $B_\lambda(x')$. We define $Z\define B_\lambda(x')\cap A$. Let $S''$ be a Weyl chamber in $A$ such that $\partial S''$ is opposite $\partial S$ and such that $Z$ is contained in $S''$. We refer to the basepoint of $S''$ by $z$  and write  $S'$ for the Weyl chamber at $z$ parallel to $S$.

By \ref{Prop_retraction} the retraction $r$ does not increase distances. Therefore $d(x',p)\geq d(x',r(p))$ for all $p\in\seg(x,y)$ and $r(\seg(x,y))$ is contained in $Z$.

Let $\widetilde{S}'$ be the unique $r(p)$-based Weyl chamber parallel to $S'$ and let $\widetilde{S}''$ be the unique $q$-based Weyl chamber parallel to $S''$. 
Applying  \ref{Lem_oppSectors} to the Weyl chambers $\widetilde{S}', \widetilde{S}''$ we obtain that $z\in\seg(q,r(p))$. Thus using the  triangle inequality we have that 
$$
d(q,r(p))=d(q,z)+d(z,r(p)) = d(q,z) + d(z,p)\geq d(q,p).
$$
On the other hand, since the retraction $r$ is distance diminishing, we have that $d(q,p)\geq d(q,r(p))$. This implies that $z$ is contained in the segment of $p$ and $q$.
By \ref{Prop_SegInApp} the points $p,q$ and $z$ are therefore contained in a common apartment $B$. The germ of $S'$ at $z$ is contained in the segment of $p$ and $q$. Thus there exists a unique Weyl chamber in $B$ based at $z$ and which is opposite $S'$ at $z$. Observe that this Weyl chamber contains $p$.
\end{proof}

\begin{prop}\label{Prop_FC}
Let $x$ and $y$ be points in $X$. For any $c\in \binfinity X$ the following is true:
\begin{itemize}[label={(FC')}, leftmargin=*]
  \item[$\mathrm{(FC)}$] The segment $\seg(x,y)$ of $x$ and $y$ is contained in a finite union of apartments containing $c$ at infinity.
\end{itemize}
\end{prop}
\begin{proof}
Lemma \ref{Lem_tec40} combined with Proposition \ref{Prop_SegInApp} implies that there exists a point $z$ such that for each $p\in \seg(x,y)$ there is a Weyl chamber $S$ based at $z$ containing $p$ which is in $\Delta_zX$ opposite the unique representative of $c$ based at $z$. By Corollary \ref{Cor_CO}
these two are contained in a common apartment. Together with $(FC')$ the assertion follows.
\end{proof}

\begin{remark}\label{Rem_FC}
Lemma 7.4.21 in \cite{BruhatTits} says that all generalized affine buildings with $\Lambda=\R$ satisfy (FC). However it is not possible to prove (FC) for generalized affine buildings using the method of \cite{BruhatTits}, since compactness arguments of $\R$-metric spaces play a major role there.
Note that the $Y$-condition of $\Lambda$-trees implies that $\Lambda$-trees satisfy (FC).
\end{remark}

%
%
%
%
%

\section{Convexity and simplicial buildings}\label{Sec_ParkinsonRam}

After we published \cite{Convexity} and communicated it to Parkinson and Ram they independently published the preprint~\cite{ParkinsonRam} in which they give a geometrical proof of one of the implications of Proposition~\ref{Prop_existenceLSgalleries} - a key in the proof of the simplicial convexity theorem. The main idea of their proof can be adapted to thick generalized affine buildings, see \ref{Sec_convexityThm}. 
For the convenience of the reader let us illustrate the main idea of their proof.

Let $X$ denote a simplicial affine building, as for example defined on p. 76 in \cite{Brown}.
Retractions $r_{A,c}$ fixed at a germ of a Weyl chamber (which is in the simplicial setting nothing but an alcove $c$) and retractions $\rho_{A,\partial C}$ fixed at a chamber at infinity can be defined as we did in Definitions \ref{Def_vertexRetraction} and \ref{Def_retractionInfty}. One can find definitions on pages 85 and 170 of \cite{Brown} as well.

Assume $X$ to be thick, i.e. each panel is contained in at least three alcoves, and let $A$ be an apartment of $X$. Fix a special vertex $0$ in $A$ and identify the spherical Weyl group $\sW$ with the stabilizer of $0$ in $\aW$. 
Let $\Cf$ denote the fundamental Weyl chamber with respect to a fixed basis $B$ of the underlying root system $\RS$. The Weyl chamber opposite $\Cf$ in $A$ is denoted by $\Cfm$ and $\cf$ denotes the fundamental alcove.
To simplify notation, write $r$ instead of $r_{A,\cf}$ and $\rho$ for $\rho_{A,(\Cfm)}$.
The main result of \cite{Convexity} is the following theorem.

\begin{thm}\label{Thm_convexity}
With notation as above let $x$ be a special vertex in $A$. 
Then
$$
\rho(r^{-1}(\sW.x)) =\dconv(\sW.x)\cap (x+\QQ) =: A^\QQ(x).
$$
where $\QQ=\QQ(\RS^\vee)$ is the co-weight lattice of $\RS$.
\end{thm}

Convexity is defined as in \ref{Def_convex}. For details compare also \cite{Convexity}. In the following a \emph{gallery} is a sequence  $\gamma= (u, c_0, d_1, c_1, \ldots , d_n, c_n, v)$ of vertices $u$ and $v$, chambers $c_i$ and panels $d_i$ where $d_i$ is contained in $c_i$ and $c_{i-1}$. The vertex $u$ is contained in $c_0$ and $v$ in $c_n$. 
We sometimes refer to $v$ as the \emph{target} of $\gamma$. 

To prove \ref{Thm_convexity} one extends $r$ and $\rho$ to galleries and uses them to describe how the building is folded onto the fixed apartment. 

The set of all minimal galleries in $X$ of fixed type with source $0$ will be denoted by $\Ghat_t$. Denote by $\G_t$ \index{$\G_t$} the set of targets of galleries $\gamma$ in $\Ghat_t$. Notice that the elements of $\G_t$ are all vertices of the same type.
Let $x$ be a special vertex in $A$ and let $\gamma:0\rightsquigarrow x$ be a minimal gallery of fixed type $t$. Then
$$
r^{-1}(\sW.x) = \G_t.
$$
If $K\define\mathrm{Stab_{\mathrm{Aut(X)}}(0)}$ acts transitively on the set of all apartments containing $0$ 
then 
$$ r^{-1}(\sW.x) = K.x .$$

Let $A$ be an apartment, $H$ a hyperplane, $d$ an alcove and $S$ a Weyl chamber in $A$. We say that \emph{$H$ separates $d$ and $S$} if there exists a representative $S'$ of $\partial S$ in $A$ such that $S'$ and $d$ are contained in different half-apartments determined by $H$.

\begin{definition}
A gallery $\gamma= (u, c_0, d_1, c_1, \ldots , d_n, c_n, v)$  is \emph{positively folded at $i$} if 
$c_i=c_{i-1}$ and 
the hyperplane $H= \Span({d_i})$ separates $c_i$ and $\Cfm$.
A gallery $\gamma$ is \emph{positively folded} if it is positively folded at $i$ whenever $c_{i-1}=c_i$.
\end{definition}

\begin{lemma}\label{Lem_folding}
Let $X$, $A$ and $\rho:X\mapsto A$ be as above and let $\hat{\rho}$ be the extension of $\rho$ to galleries. Then the image of a gallery $\gamma \in\Ghat_t$ under $\hat{\rho}$ is a positively folded gallery in $A$ of the same type.
\end{lemma}

Hence we need to understand positively folded galleries.
Let $\hat{r}$ and $\hat{\rho}$ denote the extensions of the retractions $r$ and $\rho$ to galleries. The proof of Theorem~\ref{Thm_convexity} can be reduced to the following two propositions.

\begin{prop}\label{Prop_existenceLSgalleries}
Let $A$ be a Euclidean Coxeter complex with origin $0$ and fundamental Weyl chamber $\Cf$. Let $x$ be a special vertex in $A$ and denote by $x^+$ the unique element of $\sW.x$ contained in $\Cf$. Let $t$ be the type of a fixed minimal gallery $\gamma: 0 \rightsquigarrow x^+$.
All vertices in the convex set $A^\QQ(x)\define\{y\in \dconv(\sW.x)\cap (x+\QQ(\RS^\vee))\}$ are targets of positively folded galleries having type $t$. Conversely the target of any positively folded gallery of type $t$ with source $0$ is contained in $A^\QQ(x)$.
\end{prop}

\begin{prop}\label{Prop_existencePreimage}
Let $A$ be a fixed apartment of an affine building $X$. Fix an origin $0$ and fundamental Weyl chamber $\Cf$ in $A$. If $\gamma\subset A$ is a positively folded gallery with source $0$ of type $t$ then there exists a minimal gallery $\widetilde{\gamma}\subset X$ with source $0$ such that $\hat{\rho}(\widetilde{\gamma})=\gamma$.
\end{prop}

The proof of \ref{Prop_existencePreimage} is constructive.
Proposition~\ref{Prop_existenceLSgalleries} is a purely combinatorial property of $LS$-galleries. However the proof in \cite{Convexity} made use of a character formula for highest weight representations proven by Gaussent and Littelmann in \cite{GaussentLittelmann}. We use the remainder of this section to describe the idea of a geometric proof of the first implication of  Proposition~\ref{Prop_existenceLSgalleries}.
Even though \ref{Prop_existenceLSgalleries} is stated in terms of galleries we will now give the equivalent statement using paths and the root operators defined in \cite{LittelmannPaths}.
This is the language that also applies to $\R$-buildings as defined in \cite{Ronan} or \cite{TitsComo}.

\begin{notation}\label{Not_convRev1}
Let $(X,\App)$ be the geometric realization of a simplicial affine building. Hence $\RS$ is crystallographic and the model space $\MS$ is isomorphic to the tiled vector space $V$ underlying the root system $\RS$. Let $B$ be a basis of $\RS$ with elements indexed by $I=\{1,2,\ldots,n\}$. Denote by $\Pi$ the set of all piecewise linear paths $\pi:[0,1]\rightarrow \MS$ such that $\pi(0)=0$. The \emph{concatenation of paths} $\pi_1$ and $\pi_2$ is denoted $\pi=\pi_1 \ast \pi_2$ and defined by
\begin{equation*}
\pi(t):=\left\lbrace 
\begin{array}{ll} 
	\pi_1(2t), & \text{ if } 0\leq t\leq 1/2 \\
	\pi_1(1) + \pi_2(2t-1), &\text{ if } 1/2\leq t\leq 1 .
\end{array}
\right.
\end{equation*}
We consider paths only up to re-parameterization, i.e. paths $\pi_1,\pi_2$ are identified if there exists a continuous, piecewise linear, surjective, nondecreasing map $\phi:[0,1]\rightarrow [0,1]$ such that $\pi_1\circ\phi=\pi_2$. 
\end{notation}

For any $\alpha\in B$ let $r_\alpha(\pi)$ be the path $t\mapsto r_\alpha(\pi(t))$.
Define a function $h_\alpha:[0,1] \rightarrow \R$ by $t\mapsto \lb \pi(t),\alpha^\vee\rb$ and let $n_\alpha$ be the critical value
\begin{equation}\label{Equ_tec1}
n_\alpha:=\min\{h_\alpha(t) : t\in [0,1]\}.
\end{equation}
If $n_\alpha \leq -1$  define $t_1$ to be the minimal value in $[0,1]$ such that $n_\alpha = h_\alpha(t_1)$ and let $t_0$ be the maximal value in $[0,t_1]$ such that $h_\alpha(t) \geq n_\alpha + 1$ for all $t\in [0,t_0]$.
Choose points $t_0=s_0 < s_1 < \ldots < s_r = t_1$
such that one of the two conditions holds
\begin{enumerate}
  \item \label{Num_1}$h_\alpha(s_{i-1})=h_\alpha(s_i)$ and $h_\alpha(t)\geq h_\alpha(s_{i-1})$ for all $t\in[s_i, s_{i-1}]$ or
  \item \label{Num_2} $h_\alpha$ is strictly decreasing on $[s_{i-1}, s_i]$ and $h_\alpha(t)\geq r_\alpha(s_{i-1})$ for $t\leq s_{i-1}$.
\end{enumerate}
Define $s_{-1}=0$ and $s_{r+1}=1$ and let $\pi_i$ be the path 
$$\pi_i(t)=\pi(s_{i-1} +t(s_i-s_{i-1})) - \pi(s_{i-1}) \text{ for all } i=0,1\ldots, r+1.$$
Then $\pi=\pi_0\ast\pi_1\ast \cdots\ast\pi_{r+1}$.

\begin{definition}\label{Def_rootOperator}
Let $\alpha$ be an element of $B$. Let $e_\alpha\pi:=0$ if $n_\alpha> -1$. Otherwise, let $\eta_i:=\pi_i$ if $h_\alpha$ behaves on $[s_{i-1}, s_i]$ as in (\ref{Num_1}) and let $\eta_i:=r_\alpha(\pi_i)$ if $h_\alpha$ is on $[s_{i-1}, s_i]$ as in (\ref{Num_2}). 
Then define the \emph{root operator} $e_\alpha$ associated to $\alpha$ by
$$e_\alpha\pi=\pi_0\ast \eta_1\ast\eta_2\ast\cdots\ast\eta_r\ast\pi_{r+1}.$$
\end{definition}

\begin{lemma}{\cite[Lemma 2.1.]{LittelmannPaths}}\label{Lem_rootOperator}
Let $\alpha$ be an element of $B$ and let $e_\alpha$ be as defined in \ref{Def_rootOperator}.
If $\pi$ is a path such that $e_\alpha\pi\neq 0$ then $(e_\alpha\pi)(1)=\pi(1)+\frac{2}{(\alpha,\alpha)}\alpha$.
\end{lemma}

\begin{notation}\label{Not_convRev2}
Let us add some notation to \ref{Not_convRev1}. Let $x$ be a special vertex in $\MS$ and assume without loss of generality that $x$ is contained in the fundamental Weyl chamber $\Cf$ determined by $B$, i.e  $\Cf=\{x\in\MS : \langle x, \alpha^\vee\rangle \geq 0 \;\;\forall \alpha\in B\}$. We extended the two types of retractions $r$ and $\rho$ to galleries. Analogously it is possible to consider them as maps on paths by defining the image of a point $\pi(t)$ to be the corresponding point in the image of an alcove containing $\pi(t)$. Again denote these extensions by $\hat{r}$ and $\hat{\rho}$.
\end{notation}

\begin{definition}
We say that a path $\pi'$ is a \emph{positive fold} of a path $\pi$ if there exists a finite sequence of simple roots $\alpha_{i_j}\in B$ such that $\pi'$ is the image of $\pi$ under the concatenation of the associated root operators $e_{\alpha_{i_j}}$.
\end{definition}

\begin{remark}
We could as well define \emph{positively folded paths} similarly to positively folded galleries using the notion of a billiard or Hecke path. These were defined by Kapovich and Millson in \cite{KapovichMillson}. They remark that a consequence of their results is that the so called \emph{Hecke paths} defined in \cite[3.27]{KapovichMillson} correspond precisely to the positively folded galleries defined in \cite{GaussentLittelmann}.

In particular, by Theorem 5.6 of \cite{KapovichMillson}, the LS-paths defined in \cite{LittelmannPaths} are a certain subclass of the Hecke paths.
Furthermore it is proven in \cite{KapovichMillson} that a path in an apartment $A$ of $X$ is a Hecke path if and only if it is the image of a geodesic segment in $X$ under a ``folding'' which is nothing else than a retraction centered at an alcove or a germ of a Weyl chamber. (Compare also Lemma 4.3 and 4.4 and Theorem 4.16 of \cite{KapovichMillson}.)

By \cite{LittelmannPaths} the set of $LS$-paths is invariant under the action of the root operators $e_\alpha$ with $\alpha\in B$, hence any image of a geodesic line (which trivially satisfies the axioms of an $LS$-path) obtained by applications of root operators is again an $LS$-path and therefore a Hecke path.
\end{remark}

 \begin{figure}[htbp]
 \begin{center}
 	\resizebox{!}{0.30\textheight}{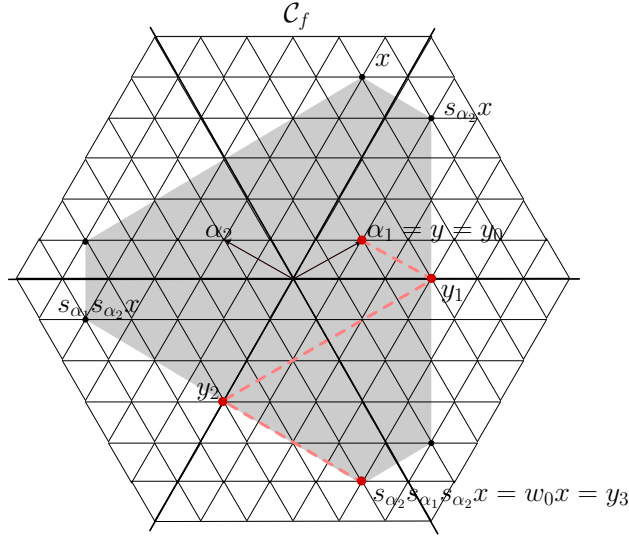}
 	\caption[adjacent]{Illustration of Lemma~\ref{Lem_ParkRam} with $x$ and $y$ as shown and  $w_0$ presented as  $s_{\alpha_2}s_{\alpha_1}s_{\alpha_2}$. The defined constants are $m_1=1, m_2=3$ and $m_3=2$.}
 	\label{Fig_folding}
 \end{center} 
 \end{figure}

\begin{lemma}{\cite[Lemma 3.1]{ParkinsonRam}}\label{Lem_ParkRam}
We denote by $\QQ^+$ the positive cone in the co-root-lattice $\QQ(\RS^\vee)$ and identify the co-root $\alpha^\vee$ of $\alpha$ with $\frac{2}{(\alpha,\alpha)}\alpha$. Let $x$ be as in \ref{Not_convRev2} and let $y$ be contained in the intersection $\bigcap_{w\in \sW } w(x - \QQ^+)$. Fix a presentation $w_0=s_{i_1}\cdots s_{i_n}$ of the longest word $w_0\in\sW$ and denote by $\alpha_{i_k}$ the root corresponding to $s_{i_k}$. Define vertices $y_i$ in the convex hull $\dconv(\sW.x)$ inductively by $y_0=y$ and for all $k=1,2,\ldots,n$ by the recursive formula 
$$
y_k = y_{k-1} - m_k\alpha^\vee_{i_k} \;\text{ where }  
$$
$$ 
  m_k=\max\{m\in\Z : y_{k-1} - m \alpha^\vee_{i_k} \in \dconv(\sW .x)\cap (x+\QQ) \}.
$$
Then $y_n=w_0x$.
\end{lemma}

Figure~\ref{Fig_folding} provides an example for the previous lemma.

In the following proposition we will restate the assertion of \ref{Prop_existenceLSgalleries} using paths instead of galleries.

\begin{prop}\label{Prop_existenceLSpaths}
Denote by $w_0$ the longest word in the spherical Weyl group. With notation as in \ref{Not_convRev2} let $\pi: 0 \rightsquigarrow x^+$ be the unique geodesic from $0$ to $w_0x$. To a special vertex $y$ in $\MS$ there exists a positive fold $\pi'$ of $\pi$ with endpoint $\pi'(1)=y$ if and only if $y$ is contained in the set
$
A^\QQ(x)=\{z\in\MS : z\in \dconv(\sW.x)\cap (x+\QQ)\}.
$
\end{prop}

The main idea of the combinatorial proof of Proposition \ref{Prop_existenceLSpaths} is to reverse Lemma~\ref{Lem_ParkRam} and ``shift'' a minimal path  (respectively gallery) $\pi$ from $0$ to $w_0x$ to a folded path from $0$ to $y$. 
Let $\pi_n\define \pi$. We apply reverse induction on $k$ ranging from $n$ to $1$. In step $k$ we define the path $\pi_{k-1}$ by an $m_k$-fold application of the root operator $e_{\alpha_{i_k}}$ to the previous path $\pi_{k}$.
According to Proposition~\ref{Lem_rootOperator} the endpoint of $\pi_k$ is translated by $\frac{2}{(\alpha_{i_k},\alpha_{i_k})} m_k\alpha_{i_k}$ in step $k$. We won't give details but illustrate the idea in the following example.

\begin{figure}[h]
\begin{center}
  \begin{minipage}[b]{7 cm}
    	\resizebox{!}{0.18\textheight}{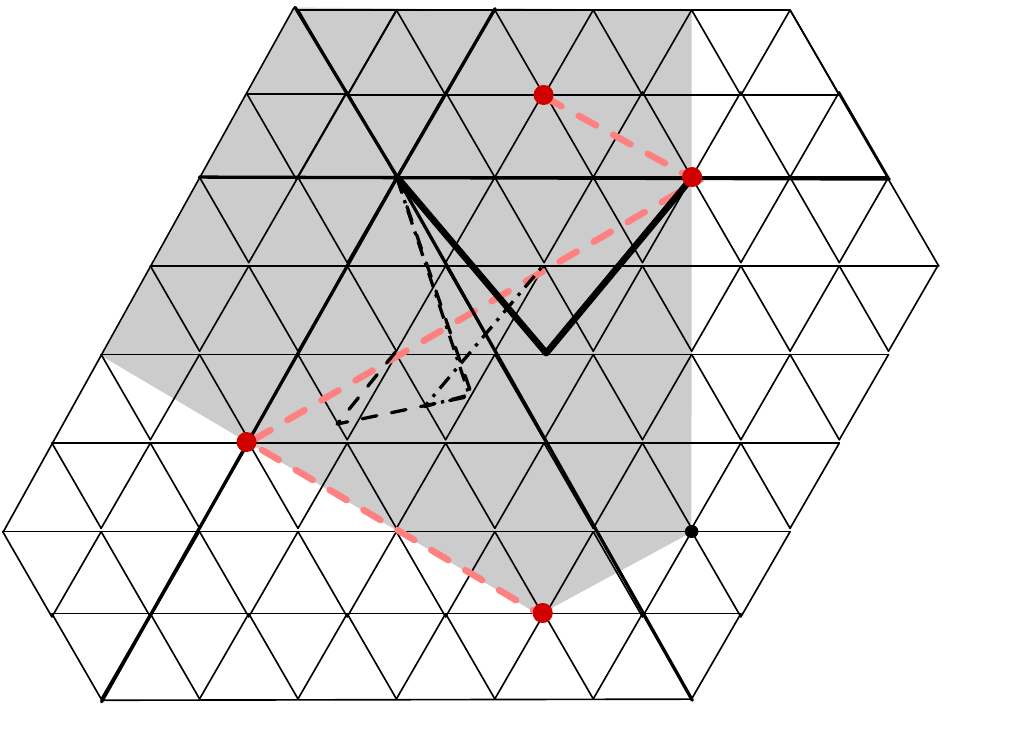}
  \end{minipage}
  \begin{minipage}[b]{7 cm}
    	\resizebox{!}{0.18\textheight}{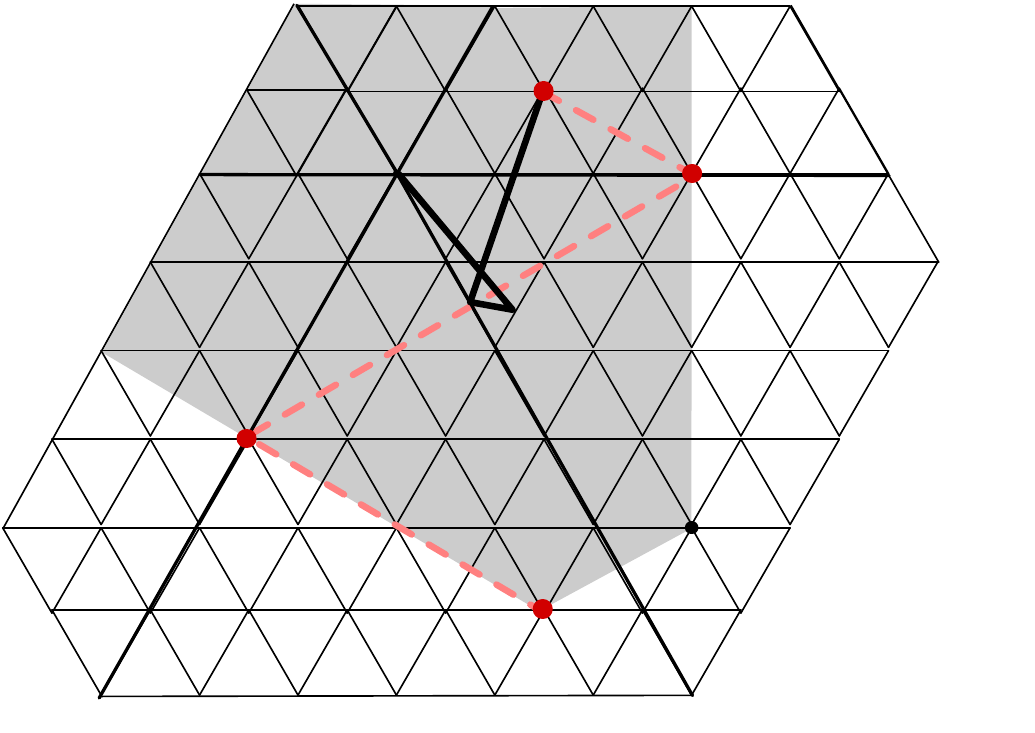}
  \end{minipage}
  \caption[adjacent]{Shifted paths according to the algorithm in the proof of Proposition \ref{Lem_rootOperator} with $w_0=s_{\alpha_2}s_{\alpha_1}s_{\alpha_2}$. The constants $m_k$ are $m_1=1, m_2=3$ and $m_3=2$.}
  \label{Fig_foldedpaths}
\end{center} 
\end{figure}

\begin{example}
In Figure \ref{Fig_foldedpaths} we illustrate the paths one obtains by applying the algorithm of the proof of \ref{Lem_rootOperator} to the example given in Figure \ref{Fig_folding}.

The vertex $x$ was chosen in $\Cf$ and an element $y\in\dconv(\sW.x)$ is given. We construct a positive fold of the minimal path connecting $0$ and $w_0x$ having endpoint $y$.

Picture (1) of Figure~\ref{Fig_foldedpaths} shows the minimal path $\pi=\pi_3$ connecting $0$ and $w_0x$. According to \ref{Lem_rootOperator} we calculate parameters $m_i$ with respect to $w_0=s_{\alpha_2}s_{\alpha_1}s_{\alpha_2}$. Observe that $m_3=2$ and that we have to apply $e_{\alpha_2}$ twice in the first step. In picture (2) the path $e_{\alpha_2} \pi$ is drawn with a dashed line, the resulting path  $\pi_2$, with notation as above, equals $e^{m_3}_{\alpha_2} \pi$ and has endpoint $y_2$. It is drawn with a solid line. Both paths are positive folds of $\pi$. Continue with a $3$-fold application of $e_{\alpha_1}$ to $\pi'$, as illustrated in picture (3). One obtains a path $\pi_1$ whose endpoint is $y_1$. Finally a single application of $e_{\alpha_2}$ gives us the path pictured in (4) which is the desired positive fold of $\pi$ with endpoint $y$.
\end{example}

\begin{figure}[h]
\begin{center}
 \begin{minipage}[b]{7 cm}
    \resizebox{!}{0.18\textheight}{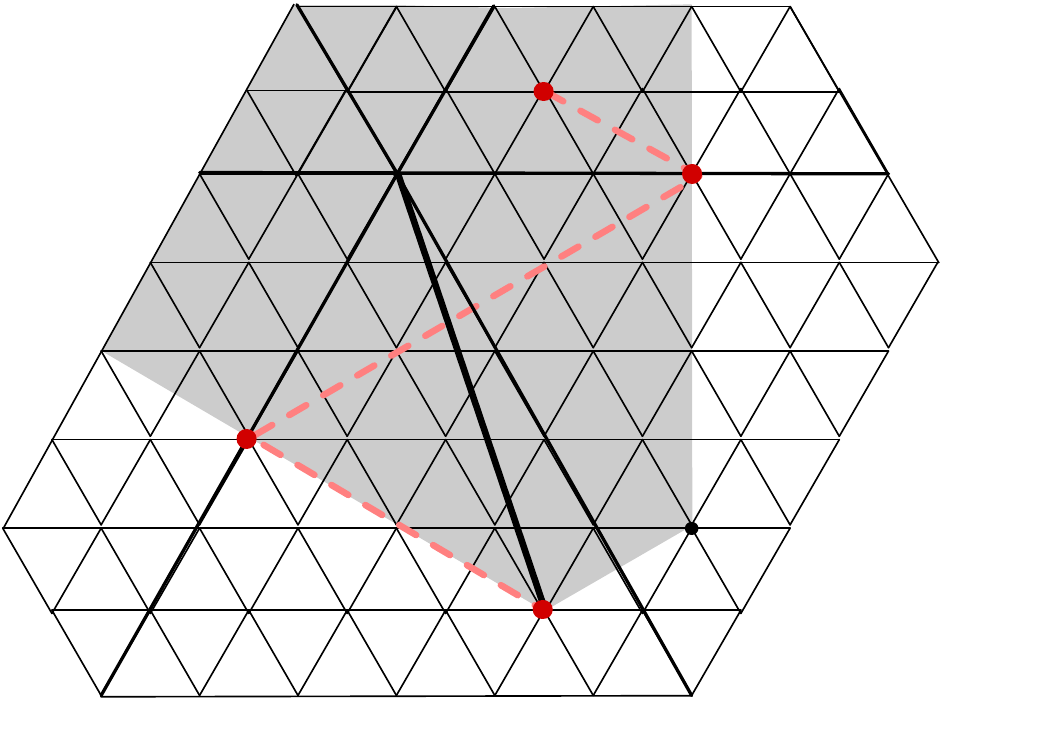}
  \end{minipage}
  \begin{minipage}[b]{7 cm}
    	\resizebox{!}{0.18\textheight}{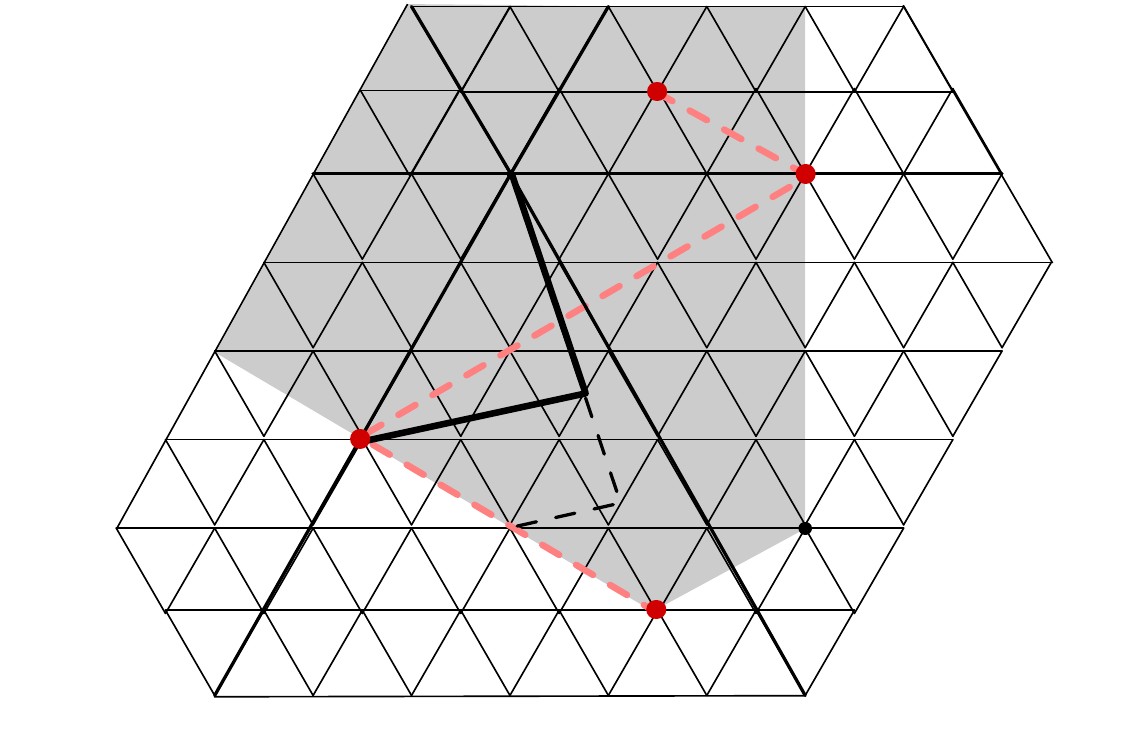}
  \end{minipage}
\end{center} 
\end{figure}

%
%
%
%

\newpage
\section{The main result}\label{Sec_convexityThm}\label{Sec_convexityRevisited}

The proof of the main result is based on the ideas discussed in Section \ref{Sec_ParkinsonRam}.

\begin{definition}\label{Def_convex}
A \emph{dual hyperplane} in the model space $\MS(\RS,\Lambda)$ is a set
$$\Hdual_{\alpha,k}=\{x\in V : ( x,\fcw ) = k\}$$
where $k\in \Lambda$ and $\fcw$ is a fundamental co-weight of $\RS$. 
Dual hyperplanes determine \emph{dual half-apartments} $\Hdual_{\alpha,k}^\pm$.
A \emph{convex set} is an intersection of finitely many dual half-apartments in $A$, where the empty intersection is defined to be $A$. The \emph{convex hull} $\dconv(C)$ of a set $C$  is the intersection of all dual half-apartments containing $C$.
\end{definition}

\begin{definition}
Fix a basis $B$ of $\RS$. The \emph{positive cone} $\Cp$ in $\MS$ with respect to $B$ is the set of all linear combinations of roots $\alpha\in B$ with coefficients in $\Lambda_{\geq 0}$.
\end{definition}

\begin{lemma}\label{Lem_conv^*Lambda}
For any special vertex $x$ in $\MS(\RS,\Lambda,T)$ let $x^+$ denote the unique element of the orbit $\sW.x$ which is contained in $\Cf$. Then
$$
\dconv(\sW.x)\cap T = \{y\in V :  x^+-y^+\in (\Cp\cap T) \}=\bigcap_{w\in \sW } w(x^+ - (\Cp\cap T)) .
$$ 
\end{lemma}
\begin{proof}
We may assume without loss of generality that $x\in \Cf$. Abbreviate $\dconv(\sW.x)\cap T$ by $A^{T}(x)$. 
Since $A^{T}(x)$ is $\sW$-invariant it suffices to prove that
$$
A^{T}(x)\cap \Cf = \{y\in \Cf : x-y \in (\Cp\cap{T}) \}.
$$
For all $\alpha\in B$ define $k_\alpha$ implicitly by $x \in \Hdual_{\alpha,k_\alpha}$.
By definition 
\begin{equation}\label{Equ_no1}
A^{T}(x) \cap\Cf = \bigcap_{\alpha\in B} (\Hdual_{\alpha,k_\alpha})^- \cap \Cf.
\end{equation}

We assume first that $y\in A^{T}(x)\cap\Cf$.  Using (\ref{Equ_no1}) and the fact that the fundamental Weyl chamber $\Cf$ is contained in the positive cone $\Cp$ we conclude that
$$
0\leq ( \fcw , y )\leq k_\alpha \text{ for all } \alpha\in B.
$$
Therefore $( \fcw,x-y ) \frac{1}{(\alpha,\alpha)} \geq 0$ and $x-y\in {T}$ since $A^{T}(x)\subset (x+{T})$. 
The faces of $\Cp$ are contained in dual hyperplanes of the form $\Hdual_{\alpha,0}$ with ${\alpha\in B}$. Therefore
$$
A^{T}(x)\cap\Cf \subset \{y\in V : x-y\in(\Cp\cap{T})\}\cap\Cf = (x-({T}\cap\Cp))\cap\Cf.
$$
Conversely assume that $y\in (x- ({T}\cap\Cp))\cap\Cf$. Then
$0\leq (\fcw, (x-y))\frac{1}{(\alpha,\alpha)}$ and hence  
$( \fcw, x ) \frac{1}{(\alpha,\alpha)} \geq  ( \fcw, y ) \frac{1}{(\alpha,\alpha)}$. Therefore the assertion holds. 
\end{proof}

The remainder of the present subsection is devoted to the proof of Theorem~\ref{Thm_convexityGeneral} which is split into two propositions \ref{Prop_preimage} and \ref{Prop_image}. For convenience let us first fix notation.

\begin{notation}\label{Not_convexity}
Let $(X,\App)$ be a thick affine building in the sense of Definitions \ref{Def_LambdaBuilding} and \ref{Def_thick} which is modeled on $\MS=\MS(\RS,\Lambda)$ equipped with the full affine Weyl group $\aW$.  Let $A$ be an apartment of $X$ and identify $A$ with $\MS(\RS,\Lambda)$ via a given chart $f$ in $\App$. Hence an origin $0$ and a fundamental Weyl chamber $\Cf$ are fixed. For simplicity we write $r$ for the retraction onto $A$ centered at $\Delta_0\Cf$ and we denote by $\rho$ the retraction onto $A$ centered at $\partial (\Cfm)$.
 As always we will identify $\sW$ with the stabilizer of the origin in $W$. 
\end{notation}

\begin{lemma}\label{Lem_key}
Let notation be as in \ref{Not_convexity} and let $x$ be a vertex in $A$. Let $w_0=s_{i_1}\cdots s_{i_n}$ be a reduced presentation of the longest word in $\sW$. Denote by $\alpha_{i_k}$ the root corresponding to $s_{i_k}$. Given a vertex $y\in \dconv(\sW.x)$ inductively define vertices $y_i$ in $\dconv(\sW.x)$ by putting $y_0=y$ and defining 
$$
y_k = y_{k-1} - \lambda_k\frac{2}{(\alpha_{i_k},\alpha_{i_k})}\alpha_{i_k}
$$ 
for all $k=1,2,\ldots,n$ with $\lambda_k\in\Lambda$ maximal such that $y_{k-1} - \lambda_k\frac{2}{(\alpha_{i_k},\alpha_{i_k})} \alpha_{i_k}$ is contained in $\dconv(\sW .x)$.
Then $y_n=w_0x^+$.
\end{lemma}
\begin{proof}
Assume without loss of generality that $x\in\Cf$. Recursively define vertices $x_0,x_1,\ldots,x_n$ in $\dconv(\sW.x)$ 
by putting $x_0=x$ and $x_k= r_{i_k}(x_{k-1}) = x_{k-1} - \lb x_{k-1},\alpha^\vee_{i_k}\rb\alpha_{i_k}$ for all $1\leq k \leq n$. 

We define a partial order on $A$ by setting $y\prec x$ if and only if $y-x\in\Cp$. 
We will show that $y_k\prec x_k$ for all $k=0,1,\ldots,n$. Then $y_n\prec x_n=w_0x$ and $y_n\in\dconv(\sW.x)$ and thus the result of the Lemma follows. We have that $y_0=y\prec x=x_0$ and by induction hypothesis that $x_{k-1} - y_{k-1}$ is contained in $\Cp$. Therefore
$$
x_k - y_k 
= x_{k-1} - y_{k-1} +(\lambda'_k - \lb  x_{k-1},\alpha_{i_k}^\vee\rb) \alpha_{i_k} 
= z + (\lambda'_k + c - \lb  x_{k-1}, \alpha_{i_k}^\vee\rb) \alpha_{i_k} 
$$ 
where $c\in\Lambda_{\geq 0}$ and $z$ is such that $x_{k-1} - y_{k-1} = z + c\alpha_{i_k}$ and that 
$$
z \in \Span_{\Lambda_{\geq 0}}(\{\alpha: \alpha\in B\setminus\{\alpha_{i_k}\}\}).
$$
We will prove that $(\lambda'_k + c - \lb x_{k-1}, \alpha_{i_k}^\vee\rb)\geq 0 $.
With $\lb  x_{k-1}, \alpha^\vee_{i_k}\rb\geq 0$ and $\lb\alpha_j,\alpha^\vee_i\rb\leq 0$ if $i\neq j$  we can conclude
$$
0\leq \lb x_{k-1},\alpha_{i_k}^\vee\rb 
= \lb y_{k-1}+c\alpha_{i_k} + z,\alpha_{i_k}^\vee\rb 
\leq \lb y_{k-1} + c\alpha_{i_k}, \alpha_{i_k}^\vee\rb.
$$
We have
$$
y_{k-1} - (\lb x_{k-1},\alpha^\vee_{i_k}\rb -c) \alpha_{i_k} 
= (y_{k-1} + c\alpha_{i_k}) - \lb x_{k-1}, \alpha_{i_k}^\vee\rb\alpha_{i_k} \in\dconv(\sW.x)
$$
and hence by definition of $\lambda_k'$ we can conclude that $\lb x_{k-1}, \alpha_{i_k}^\vee\rb -c \leq \lambda'_k$.
\end{proof}

\begin{prop}\label{Prop_preimage}
With notation as in \ref{Not_convexity} let $x$ be an element of $A$. 
For each  vertex $y$ in the convex hull $\dconv(\sW.x)$ there exists a preimage of $y$ under $\rho$ which is contained in $r^{-1}(\sW.x)$.
\end{prop}
\begin{proof}
We may assume without loss of generality that $x\in\Cf$.
We will inductively construct a preimage of $y$. Let the points $y_k$ and values $\lambda_k$ for $k=1,\ldots,n$ be as defined in the assertion of Lemma~\ref{Lem_key}.
Abbreviate $\lambda'_k:= \frac{2\lambda_k }{(\alpha_{i_k}, \alpha_{i_k})}$. The first induction step is as follows:
We denote by $H_1$ the unique hyperplane in the parallel class of $H_{\alpha_{i_1}, 0}$ containing $z_1:=y_0 - \frac{1}{2}\lambda'_1\alpha_{i_1}$. 
The associated reflection $r_1$ fixes $z_1$ and maps $y=y_0$ onto $y_1$ and vice versa.
Since $X$ is thick there exists an apartment $A_1$ in $X$ containing $y_0$ such that $A\cap A^1 = H_1^+$. We claim that $0\in H_1^+$.

The set $\dconv(\sW.x)$ is by definition $\sW$-invariant. Therefore $r_\alpha(y)$ is contained in $\dconv(\sW.x)$ for all $\alpha\in\sW$ and all $y\in\dconv(\sW.x)$. Hence $r_{\alpha_{i_1}}(y) = y- \lb y,\alpha_{i_1}^\vee \rb \alpha_{i_1}$ is contained in $\dconv(\sW.x)$ and $\lambda'_{i_1}\geq \lb y, \alpha_{i_1}^\vee \rb$ and we can conclude that $0\in H_1^+$. But $\Cf$ is contained in $H_1^+$ as well, since $\alpha_{i_1}\in B$.
Therefore $A_1$ contains $x$. We can now consider the convex hull of the orbit of $x$ under $\sW$ in the apartment $A_1$, which we denote by $\dconv_{A_1}(\sW.x)$ in order to be able to distinguish between the convex sets in the different apartments. Notice that $r$ restricted to $A_1$ is an isomorphism onto $A$. We denote the unique preimage of $y_1$ in $A_1$ by $y_1^1$ and for all $k$ the unique preimage of $y_k$ in $A_1$ by $y_k^1$.
Hence we have ``transferred'' the situation to $A_1$. 

The image of $y^1_n$ under $r$ is, by construction and Lemma~\ref{Lem_key}, the point $w_0.x$. Now the set $(A\!\setminus\! A_1 )\cup (A_1\setminus A)$ is again an apartment of $X$ containing a sub-Weyl chamber of $\Cfm$. By construction $\rho(y^1_1)=y$.

Repeating this step with $A_k$ and $y^k_k$ in place of $A$ and $y$ for all $k=1,\ldots,n$ we get a sequence of apartments $A_k$ for $k=1,\ldots,n$ with the property that they are all isomorphically mapped onto $A$ by $r$. Hence Lemma~\ref{Lem_key} implies that $y_n^n$ is a preimage of $w_0x$. Analogously to the first step observe that $\rho_{A_{k-1},\partial(\Cfm)}(y^k_k)=y^{k-1}_{k-1}$. With Lemma~\ref{Lem_compRetractions} we finally conclude that $\rho(y^n_n)=y$.
\end{proof}

\begin{property}{\bf Observation}\label{Prop_observation}
Let the notation be as in \ref{Not_convexity}. By \ref{Prop_FC'} and \ref{Prop_FC} generalized affine buildings satisfy conditions (FC') and (FC).
Let $y$ be an element of $r^{-1}(\sW.x)$. By (FC) there exists a finite collection of apartments $A_0=A$, $A_i$, $i=1,\ldots, N$ and a sequence of points $0=x_0, x_1, \ldots, x_{N-1}, x_N=y$ chosen such that 
\begin{itemize}
  \item for all $i$ the apartment $A_i$ contains a sub-Weyl chamber of $\Cfm$ and 
  \item their union contains $\seg(0,y)$ 
\end{itemize}
and further such that 
\begin{itemize}
  \item $d(0,y)=\sum_{i=0}^{N-1} d(x_i,x_{i+1})$
  \item $x_i$ and $x_{i+1}$ are contained in $A_i$ for all $i\neq N$ and such that
  \item there exist $x_i$-based Weyl chambers $S_{i-}$ and $S_{i+}$ for all $0< i< N$ containing $x_{i-1}$ and $x_{i+1}$, respectively, such that $\xi_i\define \Delta_{x_i}S_{i-}$ and $\eta_i\define\Delta_{x_i}S_{i+}$ are opposite in $\Delta_{x_i}X$. 
\end{itemize}
Let $\eta_0$ and $\xi_N$ be defined analogously.

Without loss of generality we may assume further that for all $i$ the point $x_{i-1}$ is not contained in $A_i$. Otherwise, if for some $i_0$ the point $x_{i_0}$ is contained in $A_{i_0-1}$, we define a shorter sequence of points by setting $x'_i\define x_i$ for all $i< i_0$ and $x'_i\define x_{i+1}$ for all $i\geq i_0$, that is  $x_i$ is omitted. Similarly we define apartments $A_i'$ by omitting $A_{i-1}$.

Notice that for all $i$ the restriction of  $\rho$ to $A_i$ is an isomorphism onto $A$. Hence for all $i\neq N$ the distance $d(x_i, x_{i+1})$ equals $d(\rho(x_i), \rho(x_{i+1}))$.

Let $z$ be a point contained in the interior of $\Cfm\cap \bigcap_{i=0}^N A_i$. Then by Proposition~\ref{Prop_compRetractions} $r_z\define r_{A,\Delta_z\Cf}$ equals $\rho$ on $\seg(0,y)$ and the restriction of $r_z$ to $A_i$ is an isomorphism onto $A$ for all $i=0,\ldots,N$.
\end{property}

\begin{corollary}\label{Cor_rho-distancediminishing}
The retractions $r_{A,\mu}$ and $\rho_{A,c}$ are distance non-increasing for all $A$ and $\mu\subset A$, respectively $c\in\partial A$.
\end{corollary}
\begin{proof}
Without loss of generality we may assume that $\mu=\Delta_0\Cf$. We use notation of \ref{Prop_observation}. There we made the observation that $r_{A,\mu}$ restricted to $A_i$ is an isomorphism onto $A$ for all $i$ and hence that $d(x_i, x_{i+1})=d(r_z(x_i), r_z(x_{i+1}))$ for all $i\neq N$. Therefore $d(r(x),r(y))\leq d(x,y)$ using the triangle inequality.
By Proposition \ref{Prop_compRetractions} the assertion holds for $\rho_{A,c}$ as well.
\end{proof}

\begin{lemma}\label{Lem_germs}
With notation as in \ref{Prop_observation} define $\xi'_i\define r_z(\xi_i)$ and $\eta'_i\define r_z(\eta_i)$.
Then  
$\xi'_i$ and $\eta'_i$ are chambers in $\Delta_{r_z(x_i)}A$ and there exists $w\in\sW$ such that $\xi'_i$ is the image of $w_0\eta'_i$ by $w=s_{i_1}\cdots s_{i_l}$.
\end{lemma}
\begin{proof}
By construction the germs $\xi_i$ and $\eta_i$ are opposite at $x_i$ and are contained in the apartments $A_{i-1}$ and $A_{i}$, respectively. Let $f_{i-1}, f_i$ be charts of $A_{i-1}$ and $A_i$ and assume without loss of generality that there exist a point $p_i\in\MS$ such that $f_{i-1}(p_i)=f_i(p_i)= x_i$. By axiom $(A2)$ there exists $w\in\aW$ such that 
$$
f_{i-1}\vert_{f_{i-1}^{-1}(f_i(\MS))} = (f_{i}\circ w) \vert_{f_{i-1}^{-1}(f_i(\MS))}.
$$
In fact by the assumption that $f_{i-1}(p_i)=f_i(p_i)= x_i$ the translation part of $w$ is trivial and $w$ is contained in $\sW$.
The opposite germ of $f_{i}^{-1}(\eta_i)$ is mapped onto $f_{i-1}^{-1}(\xi_i)$ by $w$.

We denote the restriction of the projection $r_z$ to $A_{i-1}$, respectively to $A_i$, by $\iota_{i-1}$ and $\iota_i$. Then $\iota_j: A_j\rightarrow A$ is an isomorphism for $j=i-1,i$. Therefore
$$
f_{i-1}^{-1}\circ \iota_{i-1}^{-1}\circ\iota_i\circ f_i : \MS\rightarrow \MS
$$
is an isomorphism, which coincides with $f_{i-1}^{-1} \circ (f_{i}\circ w)$ in a neighborhood of $x_i$.

We fix a chart $h$ of $A$ such that $h(p_i)=x'_i$. Then, since $\iota_i, \iota_{i-1}$ are isomorphisms, we conclude
$f_{i-1}^{-1}(\xi_i) = h^{-1}(\iota_{i-1}(\xi_i))$ and $f_i^{-1}(\eta_i) = h^{-1}(\iota_i(\eta_i))$.
Therefore $\xi'_i=\iota_{i-1}(\xi_i)$ is the image of $w_0 \eta'_i=\iota_{i}(\eta_i)$ by $w$, where $w_0$ is the longest word in $\sW$.
\end{proof}

\begin{lemma}\label{Lem_positiveFold}
Let notation be as in \ref{Not_convexity}. Let $x\neq 0$ be an element of $A$ and $H$ a hyperplane separating $0$ and $\Cfm$. If $y$ is a vertex in $\dconv(\sW.x)$ contained in the same half-apartment determined by $H$ which contains a sub-Weyl chamber of $\Cfm$, then the reflected image of $y$ at $H$ is contained in $\dconv(\sW.x)$.
\end{lemma}
\begin{proof}
Denote by $H'$ the hyperplane through $0$ parallel to $H$ and let $\alpha\in\RS$ be such that $H'$ is fixed by $r_{\alpha}$. The convex hull of $\sW.x$ is by definition $\sW$-invariant. Therefore the image of $y$ under the reflection at $H'$ is contained in $\dconv(\sW.x)$ and equals, by Proposition~\ref{Prop_aboveHyperplane}, $y-2y^\alpha\alpha$, since $y=m_\alpha+y^\alpha\alpha$ with $m_{\alpha}\in H'$. Similarly we can find $m\in H$ and $\lambda < y^\alpha$ such that $y=m+\lambda\alpha$. The image of $y$ under the reflection at $H$ equals $y-2\lambda\alpha$ which is obviously contained in $\dconv(\sW.x)$ as well.
\end{proof}

\begin{prop}\label{Prop_image}
Let notation be as in \ref{Not_convexity} and choose an element $x$ of $A$. The image $\rho(y)$ of $y$ in $r^{-1}(\sW.x)$  is contained in $\dconv(\sW.x)$.
\end{prop}
\begin{proof}
Let $y\in r^{-1}(\sW.x)$ be given and let $x_i, A_i$ and $S_i$ be as defined in \ref{Prop_observation}. If $N=1$ then $A_0$ contains $0=x_0, y=x_1$ and a sub-Weyl chamber of $\Cfm$ and hence $\rho(y)$ is obviously contained in $\dconv(\sW.x)$. 
By induction on $N$ and Lemma \ref{Lem_compRetractions} it is enough to prove the assertion in the case that $N=2$.
We re-use the notation of Lemma~\ref{Lem_germs} and its proof.
Recall that $A_0$ contains $x_0=0$, $x_1$ and a sub-Weyl chamber $C_0$ of $\Cfm$.

Let $(c_1=\xi_1, c_2, \ldots, c_k=\eta_1)$ be a minimal gallery in $\Delta_{x_1}X$. If $x_1$ is not contained in a bounding hyperplane of $A_{0}\cap A_1$ we may replace $x_1$ by another point satisfying this assumption.
Hence without loss of generality we can assume that $\eta_1\subset (A_1\!\setminus\! A_{0})$ and $\xi_1\subset (A_{0}\!\setminus\! A_1)$. 
Then there exists an index $j_0\in\{2,\ldots, k\}$ such that $c_{j_0-1}$ is contained in $A_{0}\!\setminus\! A_1$ and $c_{j_0}$ in $A_1\!\setminus\! A_{0}$. 

Denote by $H_1$, respectively $H_{0}$, the hyperplane spanned by the panel $p_{j_0}\define c_{j_0-1}\cap c_{j_0}$ in $A_1$, respectively $A_{0}$. Notice that $x_1$ is an element of $H_1 \cap H_{0}$.

Assume a) that \emph{$H_1$ separates $x_{2}$ and $\partial (\Cfm)$ in $A_1$ and that the point $0$ and the Weyl chamber $\partial(\Cfm)$ are separated by $H_{0}$ in $A$.}

The image of $(c_{j_0},\ldots, c_k=\eta_{j_0})$ under $\rho$ is the unique gallery $(c'_{j_0},\ldots, c'_k=\eta'_{j_0})$ of the same type which is contained in $A$ and starts in $c_{j_0-1}$. 
Hence the segment of $x_1$ and $ x_{2}$ is mapped onto $\seg_{A}(x_1, \rho(x_{2}))$ which has initial direction $\eta'_1$.
Define $x'_{2}\define\rho(x_{2})$.
By assumption the hyperplane $H_{0}$ separates $0$ and $\Cfm$. Apply Lemma \ref{Lem_germs} and obtain that  there exists $w\in\sW$ such that $\eta'_1=\rho(\eta_1)$ is the image of $w_0\xi_1$. 
Hence $x'_{2}$ is the reflected image of $r(x_{2})$ by a finite number of reflections along hyperplanes containing $x_1$ and separating $C_0$ and $0$.
Therefore $x'_{2}$ is obtained from $r(x_{2})$ by a positive fold in $A$. By Lemma \ref{Lem_positiveFold} $x'_2$ is contained in $\dconv(\sW.x)$.

In the second case b) assume that \emph{ $H_1$ separates $x_{2}$ and $\Cfm$ and assume further that  $0$ and $\Cfm$ are not separated by $H_{0}$. }

Then $\rho$ maps $(c_{j_0},\ldots, c_k=\eta_1)$ onto $(\rho(c_{j_0}),\ldots, \rho(\eta_{j_0}))$ which is a gallery of the same type contained in $A$ and $\rho(c_{j_0})$ is the unique chamber in $\Delta_{x_1}A$ sharing the panel $c_{j_0}\cap c_{j_0-1}$ with $c_{j_0-1}$. Therefore $\xi_1$ and $\eta'_1\define\rho(\eta_1)$ are opposite in $x_1$ and $\rho(x_2)=r(x_2)$.

The case c) is that \emph{ both $x_{2}$ and $0$ are not separated from $\Cfm$ by $H_1$, respectively $H_{0}$ }. 

But this cannot occur: Let $S_{1-}$ and $S_{1+}$ be the Weyl chambers based at $x_1$ having germs $\xi_1,\eta_1$ and are contained in $A$ and $A_1$, respectively. By Corollary \ref{Cor_CO} $S_{1-}$ and $S_{1+}$ are contained in a unique common apartment $B$. The span $H_B$ in $B$ of the panel $p_{j_0}=c_{j_0-1}\cap c_{j_0}$ separates the segments $\seg_B(0,x_1)$ and $\seg_B(x_1,x_2)$ and hence separates $0$ and $x_2$ which can therefore not be contained in the same half-apartment determined by $H_B$.

Finally assume d) that \emph{$H_1$ does not separate $x_{2}$ and $\Cfm$ but that $0$ and $\Cfm$ are separated by $H_{0}$.}

In this case the germ $\gamma$ of $\Cfm$ at $x_1$ is then contained in the same half-apartment of $\Delta_{x_1}$ as $\eta_1$. By the assumption that $0$ is separated from $\Cfm$ by $H_1$ there exists a minimal  gallery in $\Delta_{x_1} A$ from $\xi_1$ to $\gamma$ containing either $p_{j_0}$ or its opposite panel $q_{j_0}$ in $\Delta_{x_1}A$. If $\xi_1$ and $\gamma$ are not opposite this gallery is unique. But this contradicts the choice of $j_0$.
\end{proof}

\begin{thm}\label{Thm_convexityGeneral}
Let notation be as in \ref{Not_convexity}. Given a vertex $x$ in $A$ we can conclude
$$
\rho(r^{-1}(\sW.x)) =\dconv(\sW.x).
$$
\end{thm}
\begin{proof}
Combining Proposition~\ref{Prop_preimage} and \ref{Prop_image} we obtain the assertion.
\end{proof}

\begin{remark}
Note that Theorem \ref{Thm_convexity} is \emph{not} a special case of \ref{Thm_convexityGeneral}. It is however a special case of the conjecture \ref{Conj_convexityGeneral}.
\end{remark}

\section{An application to groups}\label{Sec_application}

For groups acting ``nicely enough'' on generalized affine buildings the proven convexity theorem can be reformulated into a statement about intersections of certain double cosets.

Let $(X,\App)$ be a thick affine building and $G$ a group acting transitively on $X$. Obviously there is an induced action of $G$ on the spherical building $\binfinity X$ at infinity. Fix a chart $f$ of an apartment $A=f(\MS)$, an origin $0=f(0)$ and assume the following transitivity properties 
\begin{enumerate}
  \item The stabilizer $G_A$ of $A$ in $G$ is transitive on points in $A$.
  \item The stabilizer $B\define G_c$ of the equivalence class $c\define\partial(\Cfm)$ is transitive on the set of apartments containing $c$ at infinity. 
  \item Assume further that $B$ splits as $B=UT$, where $T$ is the group of translations in $A$ and $U$ acts simply transitively on the apartments containing $c$ at infinity. 
  \item The stabilizer $K\define G_0$ of the origin $0$ is transitive on the set of apartments containing the origin $0$. 
\end{enumerate}

Then $G$ decomposes as follows  
$$
 G=BK=UTK \text{  \emph{ ``Iwasawa decomposition''}}.
$$ 

\begin{example}
If $G$ is a group admitting a root datum with non-discrete valuation, then $G$ has the properties described above. Compare Proposition 7.3.1 and Theorem 7.3.4 in \cite{BruhatTits}.
\end{example}

Note that points in $X$ are, under the assumptions made above, in one-to-one correspondence with left-cosets of $K$ in $G$. To see this let $x$ be a point in $X$. Lemma~\ref{Cor_tec17} and the fact that $U$ is transitive on the set of apartments containing $c$ at infinity imply that there exists $u\in U$ such that $(u^{-1}).x$ is contained in  $A$. Let $t\in T$ be the translation mapping $0$ to $(u^{-1}).x$. Identify the origin $0$ with $K$ then $x$ corresponds to the coset $utK$. One has :
\begin{align*}
\{ \text{ points in } X\} &\stackrel{1:1}{\longleftrightarrow} \{ \text{ left-cosets of } K \text{ in } G \}\\
 0 	& \longmapsto  K \\
X\ni x 	& \longmapsto  utK  \text{ with } u,t \text{ chosen as above.}
\end{align*}
Any point of $X$ can hence be identified with a coset $utK$ with suitably chosen $u\in U$ and $t\in T$. Note that points in $A$ correspond precisely to cosets of the form $tK$ with $t\in T$. 
Furthermore it is easy to see that $\rho$ is exactly the projection that maps $utK$ to $tK$ which is contained in $A$. 

For all points $tK$ in $A$ the set $r^{-1}(\sW.tK)$ is the same as the left-$K$-orbit of $tK$, i.e.
$$
r^{-1}(\sW.tK)= KtK.
$$ 
Further Lemma \ref{Cor_tec18} and the fact that $K$ is transitive on the set of apartments containing $0$ imply that $G$ decomposes as 
$$
G=KTK.
$$

The following theorem is therefore a direct reformulation of Theorem~\ref{Thm_convexityGeneral}.

\begin{thm}\label{Thm_BNpair}
For all $tK\in A$ we have that
$$
\rho(KtK)=\dconv(\sW.tK)
$$
or, since $\rho^{-1}(t'K)=Ut'K$, equivalently
$$
\emptyset \neq Ut'K \cap KtK \:\Longleftrightarrow\: t'K \in \dconv(\sW.tK).
$$
\end{thm}

%
%
%
%
%

\section{Loose ends}\label{Sec_looseEnds}

Studying generalized affine buildings which are thick with respect to a proper subgroup $\sW T$ of the full affine Weyl group $W$ we may ask again which set one obtains by taking the image under $\rho$ of the preimage $r^{-1}(\sW.x)$. These buildings only branch at special hyperplanes in the sense of Definition \ref{Def_specialHyperplane}. We conjecture that a similar answer as in the case examined in the previous chapters can be given. 

In the following let $X,\App)$ be an affine building modeled over $\MS=\MS(\RS, \Lambda, T)$ for some proper translation subgroup $T$ of $\Lambda^n$. Let $A$ be an apartment of $X$ and identify $A$ with $\MS$ via a given chart $f\in\App$. Hence an origin $0$ and a fundamental chamber $\Cf$ are fixed. We abbreviate the retraction onto $A$ centered at $\Delta_x\Cf$ by $r$ and the retraction onto $A$ centered at $\partial (\Cfm)$ by $\rho$.

\begin{conj}\label{Conj_convexityGeneral}
With notation as above assume that $X$ is thick with respect to $\WT$ (see \ref{Def_thick}). For any special vertex $x$ in $A$ we have that
$$
\rho(r^{-1}(\sW.x)) =\dconv(\sW.x)\cap (x+T).
$$
\end{conj}

Again we could restate this result in terms of groups. 
The conjecture coincides with the theorem in the classical setting if we assume $\RS$ to be crystallographic, choose $\Lambda=\R$ and let $T$ be equal to the co-root lattice $\QQ(\RS^\vee)$ of $\RS$. In this case the conjecture was proven to be true in \cite{Convexity}

The conjecture holds as well if $X$ is one-dimensional, i.e. $X$ is a $\Lambda$-tree with system of apartments, the defining root system $\RS$ is of type $A_1$ and $\MS\cong \Lambda$.

\begin{thm}\label{Thm_Conj1}
With notation as above assume that $X$ is one-dimensional and thick with respect to $\WT$. Let $x$ be a special vertex in $A$. Then
$$
\rho(r^{-1}(\sW.x)) =\dconv(\sW.x)\cap (x+T).
$$
\end{thm}
\begin{proof}
The set $\sW.x \subset A$ consists of two elements $x^+$ and $x^-$ contained in $\Cf$ and $\Cfm$, respectively. Both are at distance $d(0,x)$ from $0$. The preimage $r^{-1}(\sW.x)$ is the set of all vertices $y\in X$ such that $d(0,y)=d(0,x)$.

Let us first prove that $\rho(r^{-1}(\sW.x))$ is contained in the convex set $\dconv(\sW.x)\cap \{x+T\}$. Fix an element $y$ of $\rho(r^{-1}(\sW.x))$ which is not contained in $A$.
If $A'$ is an apartment containing $\partial(\Cfm)$ and $y$, the restriction of $\rho$ to $A'$ is an isometry onto $A$. Let $a$ denote the end of $A'$ different from $\partial(\Cfm)$ and define $z$ to be the branch point $\kappa(\partial(\Cfm) , \partial\Cf, a)$ of $A$ and $A'$. There are two cases:  Either $z$ is contained in $\Cf\setminus\{0\}$ or in $\Cfm$. In the first case $r$ and $\rho$ coincide on $A'$ and $y$ is mapped onto $x^+$. If $z\in\Cfm$ then $y\in r^{-1}(x^-)$, the distance of $x^-$ and $z$ equals $d(y,z)$ and $\rho(y)=x^- + 2d(x^-,z)\ddefine y'$, which is obviously contained in $\dconv(\sW.x)$. The point $y'$ is the reflected image of $x^-$ at (the hyperplane) $z$. Since there are no branchings other than at special hyperplanes, the branchpoint $z$ of $A$ and $A'$ is a fixed point (fixed hyperplane) of a reflection $r=t\circ r_{\alpha}$ in $\WT$. Therefore $y'\in\dconv(\sW.x)\cap \{x+T\}$.

To prove the converse let $y$ be an element of $\dconv(\sW.x)\cap \{x+T\}$. For arbitrary $v\in\MS$ let $t_v$ denote the translation of $\MS$ mapping $0$ to $v$. The unique element of $T$ mapping $x^-$ to $y$ is, in this notation, the map $t_{y-x^-}$. Let $\alpha$ be the positive rot of $\RS$ and denote by $r_\alpha$ the associated reflection of $\MS$. Apply $t_{y+x^-} \circ r_\alpha$ to $x^-$ and observe that the image is $y$. Easy calculations imply that $z\define x^- + \frac{1}{2} (y-x^-)$ is fixed by $t_{y+x^-} \circ r_\alpha$.
Therefore $z$ is a special hyperplane and, since $X$ is thick, there exists an apartment $A'$ intersecting $A$ precisely in the ray $\overrightarrow{z \partial\Cf}$. Obviously $A'$ contains $0, x^+$ and $z$, and the restriction of $r$ to $A'$ is an isomorphism onto $A$. Denote by $y'$ the preimage $r^{-1}(x^-)$ in $A'$. By construction we have $d(y',z)=d(x^-,z)$. The apartment $A''\define (A\setminus A' )\cup (A'\setminus A)$ contains $x^-,z$ and $y'$. Observe that $\rho$ restricted to $A''$ is an isomorphism of apartments mapping $y'$ onto $y$. This proves that 
$\rho(r^{-1}(\sW.x))$ contains the intersection $\dconv(\sW.x)\cap \{x+T\}$ and we are done.
\end{proof}

\phantomsection
\renewcommand{\refname}{Bibliography}
\bibliography{literaturliste}

\begin{thebibliography}{Ron89}

\bibitem[Ben90]{BennettDiss}
C.~D. Bennett.
\newblock Affine {$\Lambda$}-buildings.
\newblock {\em Dissertation, Chicago Illinois}, 106 pp., 1990.

\bibitem[Ben94]{Bennett}
C.~D. Bennett.
\newblock Affine {$\Lambda$}-buildings. {I}.
\newblock {\em Proc. London Math. Soc. (3)}, 68(3):541--576, 1994.

\bibitem[Bro89]{Brown}
K.~S. Brown.
\newblock {\em Buildings}.
\newblock Springer-Verlag, New York, 1989.

\bibitem[BT72]{BruhatTits}
F.~Bruhat and J.~Tits.
\newblock Groupes r\'eductifs sur un corps local.
\newblock {\em Inst. Hautes \'Etudes Sci. Publ. Math.}, (41):5--251, 1972.

\bibitem[BT84]{BruhatTits2}
F.~Bruhat and J.~Tits.
\newblock Groupes r\'eductifs sur un corps local. {II}. {S}ch\'emas en groupes.
  {E}xistence d'une donn\'ee radicielle valu\'ee.
\newblock {\em Inst. Hautes \'Etudes Sci. Publ. Math.}, (60):197--376, 1984.

\bibitem[GL05]{GaussentLittelmann}
S.~Gaussent and P.~Littelmann.
\newblock L{S} galleries, the path model and {MV} cycles.
\newblock {\em Duke Math. J.}, 127(1):35--88, 2005.

\bibitem[Hit08]{Convexity}
P.~Hitzelberger.
\newblock Kostant convexity for affine buildings.
\newblock {\em arXiv:math/0701094, to appear in Forum Mathematicum}, 2008.

\bibitem[KM08]{KapovichMillson}
M.~Kapovich and J.~Millson.
\newblock A path model for geodesics in {E}uclidean buildings and its
  applications to representation theory.
\newblock {\em Groups Geom. Dyn.}, 2(3):405--480, 2008.

\bibitem[Kos73]{Kostant}
B.~Kostant.
\newblock On convexity, the {W}eyl group and the {I}wasawa decomposition.
\newblock {\em Ann. Sci. \'Ecole Norm. Sup. (4)}, 6:413--455 (1974), 1973.

\bibitem[Lit95]{LittelmannPaths}
P.~Littelmann.
\newblock Paths and root operators in representation theory.
\newblock {\em Ann. of Math. (2)}, 142(3):499--525, 1995.

\bibitem[Par00]{Parreau}
A.~Parreau.
\newblock Immeubles affines: construction par les normes et \'etude des
  isom\'etries.
\newblock In {\em Crystallographic groups and their generalizations (Kortrijk,
  1999)}, volume 262 of {\em Contemp. Math.}, pages 263--302. Amer. Math. Soc.,
  Providence, RI, 2000.

\bibitem[PR08]{ParkinsonRam}
J.~Parkinson and A.~Ram.
\newblock Alcove walks, buildings, symmetric functions and representations.
\newblock {\em arXiv:0807.3602v1}, 2008.

\bibitem[Ron89]{Ronan}
M.~Ronan.
\newblock {\em Lectures on buildings}, volume~7 of {\em Perspectives in
  Mathematics}.
\newblock Academic Press Inc., Boston, MA, 1989.

\bibitem[Sch23]{Schur}
I.~Schur.
\newblock {\em Sitzungsber. Berl. Math. Ges.}, 22:9--20, 1923.

\bibitem[Tit86]{TitsComo}
J.~Tits.
\newblock Immeubles de type affine.
\newblock In {\em Buildings and the geometry of diagrams (Como, 1984)}, volume
  1181 of {\em Lecture Notes in Math.}, pages 159--190. Springer, Berlin, 1986.

\end{thebibliography}
\bibliographystyle{alpha}

\newpage

\end{document}